%% file: arxiv_mra.tex
\documentclass[12pt]{article}
\usepackage{fullpage,graphicx,amsmath,amssymb,amsfonts,verbatim,xcolor}
\usepackage[small,bf]{caption}

\usepackage[ruled]{algorithm2e} 

\SetAlFnt{\small}
\SetAlCapFnt{\small}
\SetAlCapNameFnt{\small}
\SetAlCapHSkip{0pt}
\IncMargin{-\parindent}

\usepackage[noend]{algpseudocode}

\usepackage{rotating}
\usepackage{fancyvrb}
\usepackage{subcaption}

\usepackage{url}
\usepackage{bbm}
\usepackage{tikz}
\usetikzlibrary{intersections}
\usetikzlibrary{positioning,calc,shapes,arrows}
\usetikzlibrary{positioning, arrows.meta}
\tikzstyle{block} = [rectangle, draw,
text width=7em, text centered, minimum height=5em]
\tikzstyle{line} = [draw, -latex',line width=1mm]

\usepackage{hyperref}
\hypersetup{
colorlinks   = true,    
urlcolor     = blue,    
linkcolor    = blue,    
citecolor    = red      
}

\usepackage{cleveref}
\crefname{equation}{}{}
\Crefname{equation}{}{}

\input defs.tex

\usepackage{amsthm}
\DeclareMathOperator*{\argmin}{argmin}
\DeclareMathOperator*{\argmax}{argmax}

\theoremstyle{remark}
\newtheorem{remark}{Remark}

\title{Multiple Approximate-Response Agents (MARA):\\
Fast Near-Optimal Primal Recovery for Distributed Optimization}

\author{Tetiana Parshakova\thanks{Flatiron Institute. \texttt{tparshakova@flatironinstitute.org}}
\and Yicheng Bai\thanks{Amazon. \texttt{byicheng@amazon.com}}
\and Garrett van Ryzin\thanks{Amazon. \texttt{ryzing@amazon.com}}
\and Stephen Boyd\thanks{Stanford University. \texttt{boyd@stanford.edu}}}

\date{}

\begin{document}
\maketitle

\begin{abstract}
Dual methods are useful for distributed optimization
because they allow agent-level subproblems to be solved in parallel. 
However, achieving primal feasibility
with dual methods is a challenge; it can take many iterations to find prices that recover primal feasibility, and even with optimal dual prices
primal feasibility is not guaranteed unless special conditions like strict convexity hold. 
To address this limitation, we propose a simple primal recovery method, {\em multiple approximate-response agents} (MARA),
that is able to rapidly reduce primal infeasibility,
tolerating some degree of suboptimality. The method is agnostic to how dual prices are computed, so MARA can be applied to 
enhance any dual algorithm.
Rather than returning a single primal response to each price query,
MARA requires agents to generate multiple primal responses, each of which 
has bounded suboptimality. Because these multiple responses can be computed in parallel, 
there is no increase in the wall-clock time of the underlying dual algorithm. 
MARA then constructs a convex combination of the multiple responses by minimizing the
sum of the primal and complementary slackness residuals to produce a high-quality primal solution. 
Tests of MARA using both a price localization method
and a dual subgradient method
show that it typically converges to a feasible, near-optimal solution
in a few tens of iterations. Moreover, hyperparameters of MARA can be flexibly tuned to control the
trade-off among speed, computational budget, and degree of suboptimality of the
feasible solutions.
\end{abstract}

\clearpage
\tableofcontents
\clearpage

\section{Introduction}
\label{sec:MARA}

Distributed optimization methods seek to decompose large optimization problems into a collection of subproblems (called {\em agents}) that are linked through a set of coupling constraints. Dual methods relax the coupling constraints so that the agent subproblems can be solved in parallel on distributed hardware (multiple cores and/or workers). Distributed optimization can significantly speed up the wall-clock time to find optimal solutions and greatly expand the scale of problems that can be solved.
Moreover, the collection-of-agents topology of distributed optimization reflects the structure of many real-world problems: electric power systems are collections of generation, transmission and consumption agents coupled by the circuit equations; supply chains are collections of suppliers, manufacturers, distributors and retailers linked by inventory and flow balance constraints; banking operations are collections of deposit, lending and trading divisions linked by capital budget and systemic risk constraints, etc. The combination of improved solution speed, scalability and fit to real-world problems makes distributed optimization approaches attractive in practice. 

However, one fundamental limitation of distributed optimization is achieving primal feasibility. Specifically, relying only on dual prices to induce primal feasibility requires high precision prices, and achieving high precision prices requires many dual iterations, potentially eroding the run time gains from parallelization. Worse, for important problem classes like linear models (more generally, problems that are not strictly convex), even optimal dual prices do not guarantee that the primal responses of agents will be feasible with respect to the coupling constraints. The problem of primal recovery is a major barrier to adoption of distributed optimization.

We propose a new method for primal recovery called {\em multiple approximate-response agents } (MARA)
that seeks to rapidly reduce primal residuals by 1) 
tolerating some degree of suboptimality in agent responses, and 2) leveraging parallelization to generate
multiple near-optimal primal responses from each agent. The intuition is that a large and diverse set of near-optimal agent responses can then be blended to achieve a feasible and near-optimal primal solution.
More precisely, in each iteration of the dual algorithm, we ask all agents to submit multiple, approximately-optimal responses to the current dual prices instead of just one response. (Two methods to generate these multiple responses, one based on objective value suboptimality and one based on price perturbations, are given in Section~\ref{sec-approx-oracle}.) 
Each response is an independent problem, so the multiple responses can be generated in parallel in the same amount of wall clock time as a single response. Next,  we solve a linear program to determine a convex combination of the agents' multiple responses that minimizes the sum of primal and complementary slackness residuals, producing a blended primal solution that is closer to primal feasibility. As with the multiple responses, this auxiliary linear program can be performed as an independent parallel computation that does not impact the wall clock time of the underlying dual algorithm.

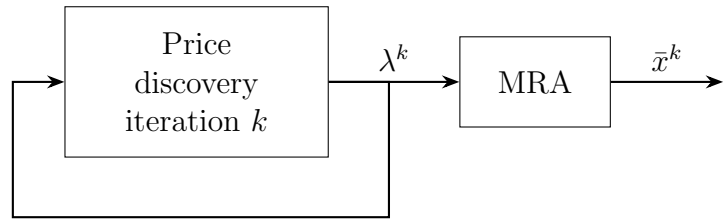
\begin{figure}[h!]
\begin{center}
\begin{tikzpicture}[auto, node distance=2cm, >=Stealth]
    \node [draw, rectangle,  align=center, minimum height=2cm, minimum width=3.5cm] 
    (B) at (0,0) {Price\\discovery\\iteration $k$};

    \node [draw, rectangle, minimum height=1.2cm, minimum width=2cm] (A) at (4.5,0) {MARA};

    \draw[->] (B.east) -- ++(0.8,0)
        [rounded corners=2pt] -- ++(0,-1.8) -- ++(-5,0) -- ++(0,1.8) -- (B.west) ;

    \draw[->] (B.east) -- (A.west) node[midway, above] {$\lambda^k$};

    \draw[->] (A.east) -- (7,0) node[midway, above] {$\bar x^k$};
\end{tikzpicture}
\end{center}
\caption{Price discovery method followed by MARA.}
    \label{fig-dual-interface-MARA}
\end{figure}

MARA can be appended to any dual algorithm, functioning as a parallel side calculation to rapidly discover feasible primal solutions while tolerating some suboptimality. 
Figure~\ref{fig-dual-interface-MARA} illustrates its logical flow:  
The first block represents
any method responsible for solving the dual problem and updating dual prices $\lambda^k$ at each iteration.
Regardless of the method used in this block,
the second block uses the current prices $\lambda^k$ and the
MARA parallel calls to agents to construct a blended primal point $\bar x^k$ with low primal residuals.

As we show in Section~\ref{sec-numerical-results}, on a range of test problems MARA is able to discover primal feasible, near-optimal solutions in tens of iterations, even though the primal solution from the underlying dual algorithm is still far from feasibility. Moreover, hyperparameters of MARA can be tuned to control the trade-off between number of iterations to find a near-optimal solution and the degree of suboptimality of that solution. 
Section~\ref{sec-MARA-extensions} discusses other variations in applying MARA that can further enhance its performance.
Our open source implementation will be available online.



The remainder of the paper is organized as follows. Section~\ref{sec-prior} provides a literature review, and the fundamentals of distributed optimization are discussed in Section~\ref{sec-distributed-optimization}. The components of MARA, including the two multiple-response oracles, are defined in Section~\ref{sec-MARA-components} and extensions and variations of MARA are presented in Section~\ref{sec-MARA-extensions}. Section~\ref{sec-numerical-results} describes our numerical tests of MARA and conclusions are given in Section~\ref{sec-conclusions}.

\section{Literature review}\label{sec-prior}

The issue of primal points not converging to feasibility in dual methods 
has long been acknowledged \cite{shor2012minimization, larsson1997lagrangean}. 
The traditional and simple solution is to perform primal averaging to 
recover primal variables in the limit \cite{larsson1999ergodic}.
For example, \cite{simonetto2016primal} use
primal averaging in a consensus-based dual decomposition framework to
obtain approximate primal solutions. 
Other works explored more sophisticated strategies.
For example, \cite{sen1986class}
combine a projected subgradient method with an auxiliary penalty
approach to produce primal-dual iterates.
And
\cite{carmon2024extracting}
propose a recursive regularization scheme
that requires an approximate primal regularized dual solve and
a query to a non-separable dual interface over the primal feasibility
set in each round
under a differentiable Lagrangian assumption. 
None of these methods, however, use the same oracles and problem
assumptions as MARA (see \S\ref{sec-distributed-optimization}) and, therefore, cannot be applied directly.

Subgradient methods for solving dual problems have been widely used in network applications to develop distributed methods \cite{low1999optimization}.
A key concern is to recover near-feasible and near-optimal primal solutions
using dual subgradients.
Primal averaging tackles this by constructing an ergodic sequence of primal variables that converges to a
primal solution asymptotically \cite{larsson1999ergodic}.
It was first proposed in the context of primal-dual subgradient schemes
\cite{nemirovskii1978cezare}
and linear problems \cite{shor1985}.
In \cite{nedic2009approximate}, the authors showed that primal averaging in
dual subgradient
methods with constant step size (as opposed to commonly
used diminishing step size)
recovers approximate primal optimal solutions.
Even though we are unaware of results applying primal averaging to dual
localization methods to recover primal solutions, 
we compare our approach to primal averaging 
in the numerical examples.

There have been decades of research on optimization methods built around
neutral cutting-plane oracles.
For example, subgradient methods were developed by Shor
in 1962~\cite{shor2012minimization}.
Extensions such as projected subgradient, Polyak's step, 
primal-dual subgradient methods, dual averaging,
and heavy ball methods have further broadened their scope;
see \cite{camerini1975improving, Polyak1987, boyd2003subgradient, 
nesterov2013introductory}.
Subgradient methods are first-order methods; thus,
their convergence is often
slow and depends on the conditioning of the problem.
However, their low per-iteration cost makes subgradient methods suitable for large-scale problems.

Another line of work that uses subgradient oracles involves localization
methods~\cite{boyd2007localization}, \ie, cutting-plane and ellipsoid methods.
These methods require more computation per iteration but typically are much
faster than subgradient methods.
Many variations have been proposed, including analytic center 
cutting-plane method~\cite{sonnevend1988new, goffin1993computation}, 
maximum volume ellipsoid cutting-plane method~\cite{tarasov1988method}, 
volumetric center cutting-plane method \cite{vaidya1989new}, and
Chebyshev center cutting-plane method~\cite{elzinga1975central}.
Using any of these methods, we can obtain optimal dual prices
conforming to our access assumptions (see \S\ref{s-price-directed-oracle}). 
Interestingly, there have been attempts to marry the two ideas of
subgradient methods with ellipsoid methods \cite{rodomanov2023subgradient},
but achieving fast practical convergence in such hybrid schemes remains an
open research question.

In our numerical experiments, we use the homogeneous analytic center
cutting-plane method (ACCPM)
\cite{nesterov1999homogeneous}.
For more details and an overview of the related methods, see \cite[\S6,7]{peton2002homogeneous}.
A cutting-plane approach typically leads to faster convergence.
This is especially beneficial when agent queries are costly,
since it allows us to quickly obtain a good estimate of the dual variable.
Specifically, we choose ACCPM because it
works well in practice and allows the approximate analytic center to be efficiently computed using the Newton method of infeasible start with diagonal Hessian
\cite{boyd2008analytic}.

\section{Distributed optimization}
\label{sec-distributed-optimization}

Consider the optimization problem 
\BEQ
\begin{array}{ll}\label{e-prob}
\mbox{minimize} & f(x)\\ 
\mbox{subject to} & Ax \leq b, 
\end{array}
\EEQ
with variable $x \in \reals^n$, 
where $f:\reals^n \to \reals \cup \{\infty\}$ is convex, closed and proper,
$A\in \reals^{m \times n}$ and $b \in \reals^m$ are given, and
the inequality is entry-wise.~\footnote{Note that affine equality
constraints can be formulated as a pair of opposing affine inequalities, 
and therefore handled in the problem~\eqref{e-prob}. Moreover, while we focus on affine constraints here for ease of exposition, more general convex constraints can be handled using MARA.}
In distributed optimization, we are interested in the special case of~\eqref{e-prob} when $f$ is block separable,
\[
f(x) = \sum_{i=1}^K f_i(x_i),  
\]
where $x = (x_1, \ldots, x_K)$,
with $f_i:\reals^{n_i} \to \reals \cup \{\infty\}$ convex, closed and
proper.
We refer to $f_i$ as the objective function of agent $i$.
Breaking the matrix $A$ into blocks conformable with $(x_1, \ldots, x_K)$,
we can then express the problem \eqref{e-prob} as
\BEQ
\begin{array}{ll}\label{e-prob-distrib}
\mbox{minimize} & \sum_{i=1}^K f_i(x_i)\\ 
\mbox{subject to} & \sum_{i=1}^K A_i x_i \leq b.
\end{array}
\EEQ

~

\noindent We refer to the separate blocks $i$ as {\em agents}. This block-separable form
makes the use of distributed optimization via dual methods appealing. 

\subsection{Price-directed oracle}
\label{s-price-directed-oracle}
We assume that $\dom f$ is bounded
and that for any $y \in \reals^n$ we can compute 
\BEQ\label{e-oracle}
x(y) \in \argmin_{z \in \dom f} \left( f(z) - y^Tz \right).
\EEQ
We refer to this oracle as \emph{price-directed}, 
since it returns an optimal action $x(y)$ given prices $y$. This oracle is typically
the primal step in a primal-dual algorithm.
When $f$ is strictly convex, the argmin is unique. However, when strict convexity does
not hold, there may be multiple optimal responses $x(y)$, in which case we assume that
the oracle \eqref{e-oracle} returns one of the optimal responses. This indeterminacy in
the oracle's response makes primal recovery challenging.
For notational convenience, we often omit writing the $y$-dependence in $x(y)$.

The price-directed oracle is closely related to $f^*$, the conjugate of $f$, as
\[
f^*(y) = \sup_{x \in \dom f} \left(y^Tx-f(x)\right) = y^Tx(y)-f(x(y)).
\]
We have $x(y) \in \partial f^*(y)$, where $\partial$ denotes subdifferential.
So, the oracle \eqref{e-oracle} is the same as finding a subgradient of
the conjugate of $f$ at $y$.  

Using the block separable form of $f$,
we can split \eqref{e-oracle}
into price-directed oracles for each $f_i$, \ie,
\BEQ\label{e-distrib-oracle}
x_i(y_i) \in \argmin_{z_i \in \dom f_i} \left( f_i(z_i) - y_i^Tz_i \right),
\quad i=1,\ldots, K,
\EEQ
where $y=(y_1, \ldots, y_K)$, conformable with $x$.
This property enables parallelization of the price-directed oracles.


\subsection{Optimality conditions}
Assuming that there exists an $x$ in the relative interior of $\dom f$ that satisfies $Ax \leq b$,
then since the problem~\eqref{e-prob} has only affine constraints,
strong duality holds from Slater's condition~\cite[\S5.2.3]{boyd2004convex}.
Therefore, the Karush-Kuhn-Tucker (KKT) optimality conditions for \eqref{e-prob} are
\BEA 
&&x = x(y), \quad y = -A^T\lambda, \label{e-opt-cond1} \\
&&Ax \leq b \label{e-opt-cond2},\\
&&\lambda_i (Ax - b)_i =0, \quad i=1, \ldots, m, \label{e-opt-cond3}\\
&&\lambda \geq 0, \label{e-opt-cond4}
\EEA
where $\lambda \in \reals^m$ is the dual variable (Lagrange multiplier) associated
with the inequality constraints $Ax \leq b$.
Assuming $\lambda \geq 0$, $y=-A^T\lambda$, and $x=x(y)$, stationarity \eqref{e-opt-cond1} and
dual feasibility \eqref{e-opt-cond4} are automatically satisfied.
Thus, the optimality conditions are reduced to
\eqref{e-opt-cond2} and \eqref{e-opt-cond3}, which are the primal feasibility and the complementary slackness, respectively.

The dual variable $\lambda$ has the interpretation of the price vector for violating 
the inequalities \eqref{e-opt-cond2}, and the variable $y=-A^T\lambda$ has the interpretation of the price vector
local to the primal variable $x$. Note that $y$ can be decomposed into price vectors
$y_i = -A_i^T \lambda$ local to each agent's primal variable $x_i$.
Define the optimal primal variable $x^\star$ and the optimal prices $\lambda^\star$.
When $f$ is strictly convex, $x^\star$ is unique, but the price vector
$\lambda^\star$ need not be.

Suppose we are given $\lambda^\star$ and $\tilde x$ satisfying
\eqref{e-opt-cond1}.
Then $\tilde x$ need not satisfy \eqref{e-opt-cond2} and \eqref{e-opt-cond3} 
unless special conditions hold, \eg, $f$ is strictly convex.
In this case, $x^\star$ is unique and therefore $\tilde x$
must be the optimal point.
However, strict convexity rules out many practically relevant cases, such as linear models.
Moreover, even in the strictly convex case, convergence can be slow because
recovering a primal solution from the dual alone often requires $\lambda^\star$ to be known
with high precision, and achieving such high precision can require a large number of dual iterations.
This behavior makes primal recovery a challenge in distributed optimization.

\subsection{Optimality condition residuals}
We define the primal and complementary slackness residuals
for a pair $(x,\lambda) \in \reals^n \times \reals_+^m$ as
\[
r_p = \ones^T (Ax-b)_+, \qquad r_c = \lambda^T |Ax - b|,
\]
where $\ones$ is the vector with all entries,
$u_+=\max\{u,0\}$ is the positive part, and the positive part
and absolute values are applied entrywise.
Evidently $r_p$ and $r_c$ are nonnegative, with $(x,\lambda)$ 
optimal (assuming $\lambda \geq 0$, $y=A^T\lambda$, and $x=x(y)$) if and only if
$r_p+r_c=0$.
We use these residuals as a stopping criterion, 
\[
r_p + r_c \leq \epsilon_r,
\]
where $\epsilon_r>0$ is a given residual tolerance.

\subsection{Dual subgradients}
The dual function associated with the problem~\eqref{e-prob} is
\[
g(\lambda)  = -f^*(- A^T\lambda) - \lambda^T b
\]
for $\lambda \geq 0$.
The optimal prices $\lambda^\star$ maximize $g$ over $\lambda \geq 0$.
The negative dual function $-g$ is convex, with subgradient
at $\lambda \geq 0$
\BEQ\label{e-dual-subgrad}
-Ax(y) + b \in \partial (-g)(\lambda),
\EEQ
which we can construct from our agent price-directed oracles. 
These subgradients are used in various primal-dual algorithms.

\section{Components of MARA}
\label{sec-MARA-components}

We next describe the two main components of MARA: 1) oracles for generating multiple responses, and 2) a primal recovery method that forms a primal solution by blending the multiple responses.

\subsection{Multiple approximate-response oracles}\label{sec-approx-oracle}

For any $y \in \reals^n$,
an {\em approximate price-directed oracle}
returns $x^{\text{apx}}(y)$ that approximately
minimizes 
$f(z) - y^{T} z$ with respect to $z$,  
\ie, 
\BEQ\label{e-appx-oracle}
    -f^*(y)  \leq f(x^{\text{apx}}(y)) - y^{T} x^{\text{apx}}(y) 
    \approx -f^*(y).
\EEQ
We next construct two approximate oracles that ensure that
any convex combination of $\epsilon$-suboptimal points
results in only a small violation
of the stationarity condition \eqref{e-opt-cond1}.

\subsubsection{Value suboptimality oracle}
For any $y \in \reals^n$, 
we define an $\epsilon_v$-value suboptimal primal variable 
with respect to \eqref{e-appx-oracle}
as $x^{\text{v}}(y)$ such that
\[
    - f^*(y) \leq  f(x^{\text{v}}(y)) - y^Tx^{\text{v}}(y)  \leq - f^*(y) +\epsilon_v |f^*(y)|.
\]

To implement this oracle, we follow the approach described in \cite{skaf2010techniques}.
Specifically, we sample $N$ i.i.d. directions
$\delta_1, \ldots, \delta_N \sim \mathcal{N}(0, I_n)$.
For each $j=1, \ldots, N$, we compute $z_j$ by solving a convex problem
\[
\begin{array}{ll}
\mbox{maximize} & \delta_j^T z_j\\ 
\mbox{subject to} & f(z_j) - y^T z_j \leq - f^*(y) +\epsilon_v |f^*(y)|.
\end{array}
\]

~

The solutions $\{z_1, \ldots, z_N \}$ form
a list of $N$ $\epsilon_v$-suboptimal
primal variables associated with local price $y$.
Further for any convex combination of these points,
$\bar x = \sum_{j=1}^N \theta_j z_j$ with $\ones^T\theta=1$ and $\theta \geq 0$, the suboptimality w.r.t. stationarity
condition \eqref{e-opt-cond1} is preserved
\[
    f(\bar x) - y^T \bar x 
    \leq \sum_{j=1}^N \theta_j(f(z_j) - y^Tz_j)
    \leq - f^*(y) +\epsilon_v |f^*(y)|.
 \]
 Note this oracle requires a modification of the agent subproblems.

\subsection{Price perturbation oracle}
For any $y \in \reals^n$, 
we also define the primal point $x^{\text{p}}(y)$ with $\epsilon_p$-perturbed prices
with respect to \eqref{e-appx-oracle}
as 
\[
    f(x^{\text{p}}(y)) - (y + \delta)^Tx^{\text{p}}(y) =- f^*(y + \delta),
\]
with perturbation vector $\delta \in \left [-\epsilon_p|y|, \epsilon_p|y| \right ]$.
To compute this, we sample $\delta_1, \ldots, \delta_N$ uniformly at random from the
bounding box
$\left [-\epsilon_p|y|, \epsilon_p|y| \right ]$.
For each $j=1, \ldots, N$, we compute $z_j \in \partial f^*(y + \delta_j)$. 
The solutions $\{z_1, \ldots, z_N \}$ form 
a list of $N$ 
$\epsilon_p$-perturbed prices
with respect to $y$.

If we assume that $f^*$ is $L$-smooth,
we can get the bound on the suboptimality of this oracle.
Recall that $f^*$ is $L$-smooth if and only if $f$ is $(1/L)$-strongly convex.
Then for any convex combination
$\bar x=\sum_{j=1}^N \theta_j z_j$ we have
\BEAS
    f(\bar x) - y^T \bar x
    &\leq& 
    \sum_{j=1}^N \theta_j(f(z_j) - y^Tz_j) \\
    &\leq& 
    \sum_{j=1}^N \theta_j(f(x) - (y + \delta_j)^T(x - z_j) -  y^Tz_j)\\
    &=& -f^*(y) + \sum_{j=1}^N \theta_j \| \delta_j\|_2 \|x - z_j\|_2 \\
    &\leq& -f^*(y) + \epsilon_p^2 L\|y\|_2^2.
\EEAS
In the above, we used Jensen's inequality for convex $f$,
subgradient inequality 
$f(z_j) \leq f(x) - (y + \delta_j)^T(x - z_j)$,
Cauchy-Schwarz inequality, $\|\delta_j\|_2 \leq \epsilon_p\|y\|_2$.
Also by $L$-smoothness of $f^*$, we have 
$\|x - z_j\|_2 = \|  \nabla f^*(y) - \nabla f^*(y+\delta_j)\|_2\leq L \|\delta_j\|_2$.
The point $\bar x$ is suboptimal with an additive error that is proportional to 
the square of the price perturbation's norm. Note that because it is price based, this oracle does not require any modification to the agent subproblems.


\subsection{Primal recovery method}\label{sec-cvx-primal-recovery}
We next define a primal recovery method that uses the multiple agent responses to 
generate a new primal variable which reduces the sum of residuals, $r_p + r_c$.
This recovery procedure is an independent post-processing step that can be applied flexibly, 
either at every iteration of the dual algorithm, applied periodically (\eg, every $10$ iterations), 
or only applied after some measure of dual prices have stabilized.

Starting from a given dual variable $\lambda \in \reals_+^m$, it consists of two steps:

\paragraph{Step 1. Approximate oracle queries}
Using an approximate oracle from \S\ref{sec-approx-oracle}, 
each agent $i=1,\ldots, K$ returns a list of $N_i$ $\epsilon$-suboptimal
primal vectors,
$\{z_{i_1}, \ldots, z_{i_{N_i}}\}$, 
associated with local price vector $y_i = - A_i^T\lambda$.
Note that this step can be computed by making parallel calls to agents,
so the wall-clock time is equivalent to a single agent response.

\paragraph{Step 2. Primal recovery}
Construct a convex combination of 
the responses of each agent, $\bar x_i$,
minimizing the sum of residuals, $r_p + r_c$, by solving the following linear programming
problem: 
\BEQ\label{e-cvx-primal-recov}
\begin{array}{ll}
\mbox{minimize} & r_p + r_c  \\
\mbox{subject to} & r_p = \ones^T (A\bar x-b)_+ \\
& r_c = \lambda^T |A\bar x - b| \\
& \bar x_i = Z_i u_i, \quad i=1, \ldots, K \\
& \ones^T u_i = 1, \quad i=1, \ldots, K \\
& u_i \geq 0, \quad i=1, \ldots, K \\
& \bar x = (\bar x_1, \ldots, \bar x_K)
\end{array}
\EEQ
with variables $u_i \in \reals^{N_i}$ for all $i=1, \ldots, K$, and
matrices
\BEQ \label{e-zi-matrices}
Z_i = \left[ \begin{array}{ccc} 
z_{i_1}& \cdots & z_{i_{N_i}}
\end{array}\right] \in \reals^{n \times N_i},
\quad i=1, \ldots, K.
\EEQ
Using results in \S\ref{sec-approx-oracle},
the primal point $\bar x$ introduces small violations
of the stationarity condition \eqref{e-opt-cond1}
while simultaneously minimizing the sum of residuals, $r_p + r_c$.

\begin{remark}
Note that since $\bar x$ is constructed from primal points obtained using
an approximate oracle,
$r_p+r_c=0$ does not guarantee that $\bar x$ is
primal optimal.
This is because $\bar x$ need not respect the stationarity condition
\eqref{e-opt-cond1}.
\end{remark}

\section{Extensions and variations}
\label{sec-MARA-extensions}
There is significant flexibility in applying MARA.
In this section, we highlight a few key variations of the method.

\subsection{Minimizing $r_p$}
First, since we are ultimately interested in reducing primal residuals,
problem \eqref{e-cvx-primal-recov} can be replaced by the one below
\BEQ\label{e-cvx-primal-recov-rp}
\begin{array}{ll}
\mbox{minimize} & r_p  \\
\mbox{subject to} & r_p = \ones^T (A\bar x-b)_+ \\
& \bar x_i = Z_i u_i, \quad i=1, \ldots, K \\
& \ones^T u_i = 1, \quad i=1, \ldots, K \\
& u_i \geq 0, \quad i=1, \ldots, K \\
& \bar x = (\bar x_1, \ldots, \bar x_K),
\end{array}
\EEQ

~

\noindent where we removed the complementary slackness residual $r_c$ from the objective in \eqref{e-cvx-primal-recov}.
Originally, the complementary slackness residual $r_c$ was included
in the objective to prevent us from generating a primal point
that is too suboptimal. 
However, due to the convexity of $f_i$ for each $i=1, \ldots, K$, if each candidate point in $Z_i$ has tolerable suboptimality with respect to $f_i(x_i^\star)$,
then so does any convex combination. 
Consequently, the problem \eqref{e-cvx-primal-recov-rp}  
will find a point $\bar x$ that also has tolerable suboptimality,
with smaller primal residuals than the point returned by \eqref{e-cvx-primal-recov}.

\subsection{Increasing the primal set with history}
If the approximate oracle generates primal points with comparable
function values in the last $H_i$ iterations, then we can add
these additional candidate points, $Z_i^{k-H_i+1}, \ldots, Z_i^k$ for all $i=1, \ldots, K$, into the candidate set of responses used to construct a new point
$\bar x^k$ using \eqref{e-cvx-primal-recov-rp}. 
Due to convexity of the agents' functions, 
the resulting point will have a comparable function value,
while searching over the larger polytope to reduce the primal residuals.

\subsection{Oracle suboptimality mixing}
The suboptimality in the approximate oracle controls the violation of the
stationarity condition \eqref{e-opt-cond1}.
In particular, let $\bar x_1$ and $\bar x_2$ be returned
by MARA primal recovery with additive suboptimalities $\epsilon_1 \gg \epsilon_2$, respectively,
\ie,
\[
-f^*(y) \leq f(\bar x_j) - y^T \bar x_j \leq -f^*(y) + \epsilon_j,
\quad j=1, 2.
\]
This implies
\[
f(\bar x_1) \gg f(\bar x_2) - \|y\|_2\|\bar x_2 - \bar x_1\|_2. 
\]
Therefore, if $\|y\|_2\|\bar x_2 - \bar x_1\|_2$ is sufficiently small,
the primal point $\bar x_2$
will have a lower function value (\ie, $f(\bar x_1) > f(\bar x_2) $),
but potentially higher residuals $r_p+r_c$.
This is because the set of $\epsilon_2$-suboptimal primal points
is a subset of the $\epsilon_1$-suboptimal points.
This suggests splitting the budget of responses
$N_i$ across
multiple approximate oracles with varying levels of suboptimality.
This approach balances between suboptimality and primal feasibility, allowing us
to quickly obtain the primal feasible point with large $\epsilon$,
while still allowing further improvement in suboptimality and
feasibility as the dual iterate improves.

\section{Numerical results}
\label{sec-numerical-results}

To test MARA we ran numerical experiments using two different dual methods applied to four problem classes. We first discuss the two dual methods, then describe the experimental process and finally discuss the results of MARA on each of the four problem classes.

\subsection{Price discovery methods}\label{sec-price-discovery}

\subsubsection{Projected dual subgradient method}

Consider first a subgradient method to solve
the dual problem \cite{boyd2003subgradient}.
Let $\mathcal{P}$ be the initial set containing optimal dual variable, 
\eg, $\mathcal{P} = \reals^m_+$.
Let $\Pi_\mathcal{P}:\reals^m \to \mathcal{P}$ denote the Euclidean 
projection onto this set.
Then the projected dual subgradient method is given by Algorithm~\ref{alg-dual-subgrad}.
Note that this algorithm requires the feasible set $\mathcal{P}$,
initial point $\lambda^0$
and the step size rule to update $\alpha^k$.

\begin{algdesc}{\sc Projected dual subgradient method}\label{alg-dual-subgrad}
    {\footnotesize
    \begin{tabbing}
    {\bf given} feasible set $\mathcal{P}$, initial point $\lambda^0$. \\
    {\bf for iteration $k=0,1,\ldots$} \\*[\smallskipamount]
     Query \eqref{e-dual-subgrad} at $\lambda^k$
     to get a $g(\lambda^k)$ and $q^k \in \partial(-g)(\lambda^k)$.\\
    Set the step size $\alpha^k$. \\
    Update iterate $\lambda^{k+1} = \Pi_\mathcal{P}(\lambda^k - \alpha^k q^k)$. \\
    Compute best value $g^k_{\mathrm{best}} = \max\{g(\lambda^k), 
    g^{k-1}_{\mathrm{best}} \}$. \\
     \emph{Check stopping criterion.} Quit if 
     $g^{k}_{\mathrm{best}} - g^{k-1}_{\mathrm{best}}$ is small.
    \end{tabbing}}
\end{algdesc}

\subsubsection{Price localization method}\label{s-price-local}

Unlike dual subgradient methods, which rely on simple steps based on excess consumption,
localization methods perform substantial calculation to update dual prices at each iteration.
As a result, they tend to produce higher-quality dual prices within a small number of iterations.
Moreover, localization methods do not rely on function evaluations,
which enables them to handle more general cases when the agents do not have 
an objective and instead use a monotone operator that maps prices
to actions.

Recall \eqref{e-dual-subgrad},
then subgradient inequality for the convex function $-g$, together
with $g(\lambda^\star) \geq g(\lambda)$, gives
\[
(-Ax(y) + b)^T(\lambda^\star - \lambda)\leq 0,
\]
\ie, any $\lambda^\star$ lies in the halfspace
\BEQ\label{e-cutting-plane}
\left\{ \tilde \lambda ~{\Big |}~ (\underbrace{-Ax(y)+b}_{=:c})^T \tilde \lambda \leq \underbrace{(-Ax(y)+b)^T\lambda}_{=:d},
~ y=-A^T\lambda
\right\} = 
\left \{ \tilde \lambda \mid c^T \tilde \lambda \leq d \right \}.
\EEQ
Thus by evaluating $x(y)$ for $\lambda \geq 0$ (which is consistent
with our access assumptions) we obtain a 
cutting plane for $\lambda^\star$ that passes through $\lambda$,
\ie, a neutral cutting plane for the optimal price.

Localization methods for solving \eqref{e-prob} use the neutral
cutting plane to refine the set that contains the optimal
price $\lambda^\star$. 
They differ in how they find the ``center'' point of the current set from which
to query the cutting plane.
We use the homogeneous analytic center 
cutting-plane method \cite{nesterov1999homogeneous},
which embeds the problem into extended space, 
normalizes the iterates via the proximal term, and
offers more flexibility in terms of initial information than related methods.
Specifically, this involves embedding the dual problem
into a conic form with
augmented variable 
$z=(t, \bar \lambda)\in \reals_+ \times \reals^m$
such that the dual variable is $\lambda=\bar \lambda/t$.
This method also employs augmented self-concordant barriers, given by the sum of a 
self-concordant barrier for the cone and a proximal term.
To summarize, the localization method is described in Algorithm~\ref{alg-localization}. 

\begin{algdesc}{\sc Homogeneous analytic center
cutting-plane method}\label{alg-localization}
    {\footnotesize
    \begin{tabbing}
    {\bf given} initial $\nu$-normal barrier $F(z)$. \\*[\smallskipamount]
    Set initial function $F^0(z) = F(z) + (1/2)\|z\|_2^2$. \\
    {\bf for iteration $k=0,1,\ldots$} \\*[\smallskipamount]
     Compute center point $(t^{k+1}, \bar\lambda^{k+1})=
      \argmin_z F^k(z)$. \\
     Query \eqref{e-cutting-plane} at $\lambda^{k+1}=\bar \lambda^{k+1}/t^{k+1}$
     to get a halfspace $\left\{ \lambda \mid (c^{k+1})^T \lambda \leq d^{k+1} \right\}$.\\
    Compute normalization $n^{k+1} = (\|c^{k+1}\|^2_2 + (d^{k+1})^2)^{1/2}$. \\
    Update the function $F^{k+1}(z) = F^k(z) - \log ((td^{k+1} - (c^{k+1})^T \bar \lambda)/n^{k+1})$, where $z=(t, \bar \lambda)$. \\
     \emph{Check stopping criterion.} Quit if 
     $\| \lambda^{k+1} - \lambda^{k} \|_2$ is small.
    \end{tabbing}}
\end{algdesc}

For the initial barrier $F(z)$ we use the logarithmic barrier 
\[
F(z) = -\sum_{j=1}^l \log (td^j - (c^j)^T \bar \lambda) - \log t,
\]
which is $(l+1)$-normal barrier function.
Algorithm~\ref{alg-localization} converges to $\epsilon$ accuracy in $O(\nu/\epsilon^2)$
iterations.
For more details see \cite[\S6,7]{peton2002homogeneous} and \cite{nesterov1999homogeneous}.
In localization methods, \eg, ACCPM, we establish the convergence of 
the best-to-date function value. 
Interestingly, in homogeneous ACCPM, we can construct a sequence of points
that ensures the convergence of the function values; see Appendix \ref{sec-haccpm}.

\subsection{Numerical testing procedure}

We used the following procedure to test MARA:
For every problem class, at every iteration $k$ the given price discovery method computes a
dual vector $\lambda^k$ using agent price-directed oracles $x^k \in \partial f^*(-A^T \lambda^k)$.
We assume we have a good approximation of the range of optimal prices, 
thus we set the initial polyhedron for the localization method
within a box that spans from $1/3$ to $3$ times the estimated
maximum and minimum price values.
We also use this box as the feasible set $\mathcal{P}$ in the dual subgradient method.
For the step size rule, 
we select the value from the set 
$\alpha^k \in \{ 0.1/\sqrt{k}, 1/\sqrt{k}, 1/k, 10/k\}$
that minimizes the final primal residuals for each example.

We then applied the MARA primal recovery method to compute $\bar x^k$ with respect to 
prices $y^k=-A^T \lambda^k$ in each iteration.
We tried both the value suboptimality and price perturbation approximate oracles for each example.
We also report the suboptimality of the function value for the
projection of the primal point $x^k$ onto the primal feasibility set
\[ \Pi (x^k) = \argmin_{x: Ax \leq b} \|x - x^k\|_2 .
\]
Lastly, since primal averaging in dual subgradient methods is a common technique to recover the primal variable in the limit \cite{larsson1999ergodic},
for comparison we also report the performance of the primal average point,  
$\frac{1}{k}\sum_{j=1}^k x^j$, as another benchmark.

We consider the solution primal feasible if the relative primal infeasibility
satisfies $\|r_p\|/\|b\|_2 < 10^{-6}$.
Results for every problem class are reported in the same four panel figure (\eg, see Figure~\ref{fig-ra-subopt-vp}), where for every experiment we plot best-to-date suboptimality of the latest feasible point
in the top left panel;
any improvement is marked with a star $\star$. If a method never recovers
a feasible point, its curve will not appear in this panel.

We initially set the suboptimality threshold for approximate response to $\epsilon=10\%$ and number of approximate responses
$N_i=10$ for all $i=1, \ldots, K$. 
In \S\ref{sec-exp-tradeoffs} we also compare the performance of MARA
for $\epsilon=10\%$ and $\epsilon=1\%$
as well as
$N_i=10$ and $N_i=50$ for all $i=1, \ldots, K$.
Lastly, in \S\ref{sec-exper-history} we discuss the performance of MARA when history
is included, \ie, $H\geq 2$.

\subsection{Resource allocation problem}
As a first example we 
consider the resource allocation problem, 
\[
\begin{array}{ll}
\mbox{minimize} & \sum_{i=1}^K f_i(x_i)\\ 
\mbox{subject to} & \sum_{i=1}^K x_i \leq R,
\end{array}
\]
where $f_i:\reals^{m}_+ \to \reals \cup \{\infty\}$ is a negative utility function, 
$R \in \reals^{m}_+$ is a vector with budgets for $m$ resources, 
and $x_i \in \reals^{m}_+$ are 
variables for all $i=1, \ldots, K$. 
Our example uses negative utility functions of the form
\[
f_i(x_i)=  
\begin{cases}
    -\mathcal G \left(C_ix_i\right) & x_i \geq 0, ~ x_i^T\ones \leq 1 \\
    \infty & \text{otherwise},
\end{cases}
\]
where $\mathcal G (u) = \left( \prod_{i=1}^p u_i\right)^{1/p}$ for 
$u \in \reals_+^p$ is the geometric mean function.
The entries of $C_i$ are nonnegative, so $-f_i$ is concave and nondecreasing.

\paragraph{Problem instance}
We choose the entries of $(C_i)_{j\cdot}$ to be uniform on $[(j-1)/m,j/m]$, 
and the resource budgets $R_{m-j}$ are chosen uniformly from
$\sqrt{n}+(n/2-\sqrt{n})[(j-1)/m,j/m]$
for all $j=1, \ldots, m$,
this way the more efficient resources
are also more scarce.
We consider $m=50$ resources, $M=10$, 
and $K=100$ agents.
The utility functions use $p=5$, \ie, each is the 
geometric mean of $p$ affine functions.
The dimension of primal variable is $5000$
and of dual variable is $50$.

\begin{figure}[t]
\begin{center}
    \centering
        \includegraphics[width=0.9\textwidth]{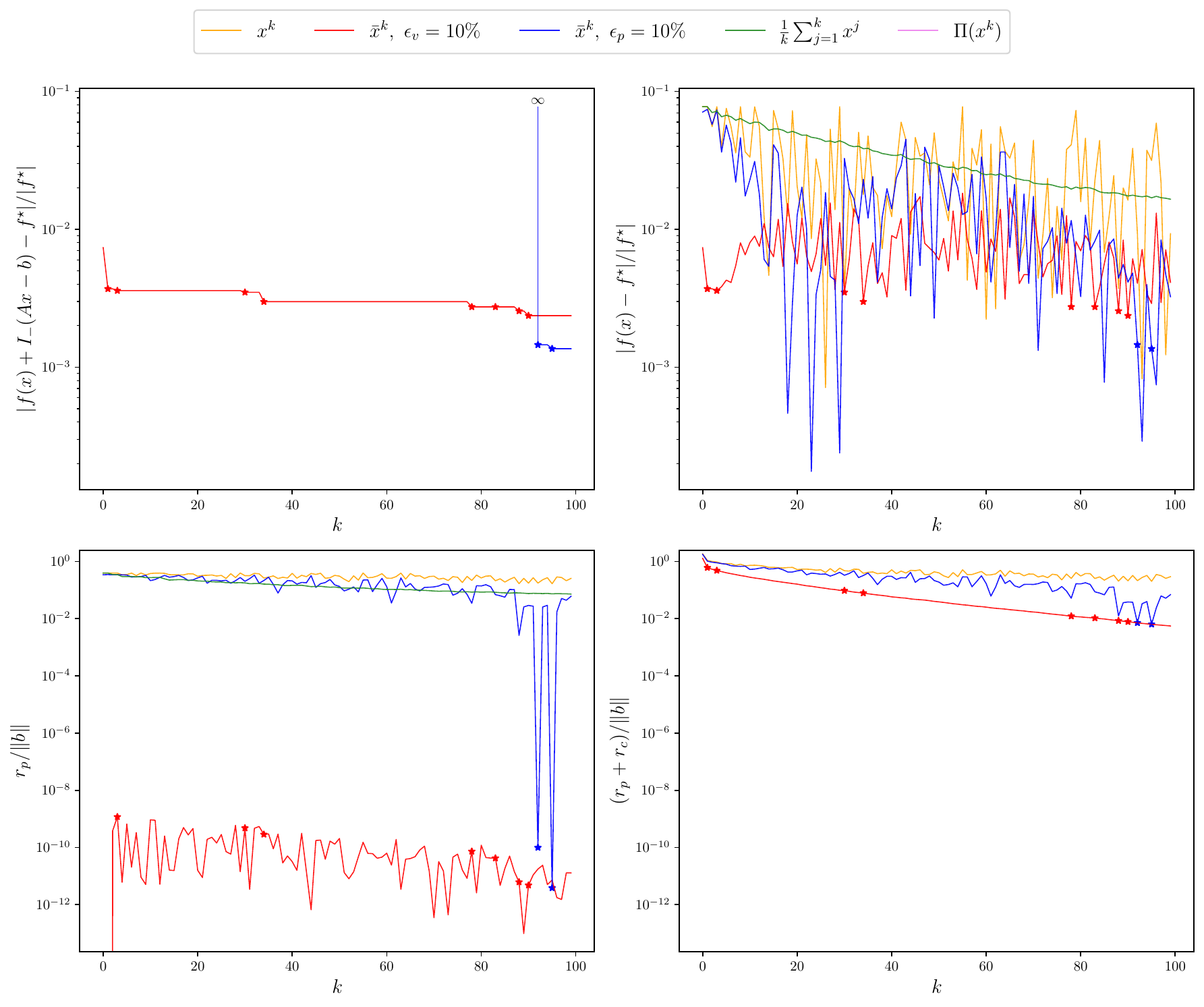}
    \caption{Best-to-date suboptimality of the most recent feasible point (top left),
    function value suboptimality (top right), relative primal violations (bottom left),
    relative residuals (bottom right)
    versus localization method iterations for the resource allocation problem.}
    \label{fig-ra-subopt-vp}
\end{center}
\end{figure}

\paragraph{Results}
Figure~\ref{fig-ra-subopt-vp} illustrates the suboptimality and relative residuals achieved with the localization method
for the resource allocation problem.
Note that in each iteration MARA using the value suboptimality oracle with $\epsilon_v=10\%$ 
recovers a feasible primal
point $\bar x^k$ (relative primal residuals in the lower left panel are in the 
$[10^{-11},10^{-9}]$ range), achieving, for example, $1.2\%$ function value suboptimality in iteration $25$.
At iteration $92$, MARA using a price perturbation oracle with $\epsilon_p=10\%$ finds a 
feasible primal
point $\bar x^k$ with a function value suboptimality of $0.1\%$.

In contrast, at iteration $95$, the original primal value $x^k$ is still primal infeasible with $16.7\%$ 
relative infeasibility
and $3.2\%$ suboptimality, and even the primal average point at $k=99$ remains infeasible with $7.2\%$
relative infeasibility
and
$1.7\%$ suboptimality.
Moreover, the projected point $\Pi(x^k)$ is everywhere outside the $\dom f$ and therefore has unbounded negative utility.

\begin{figure}[t]
\begin{center}
    \centering
        \includegraphics[width=0.9\textwidth]{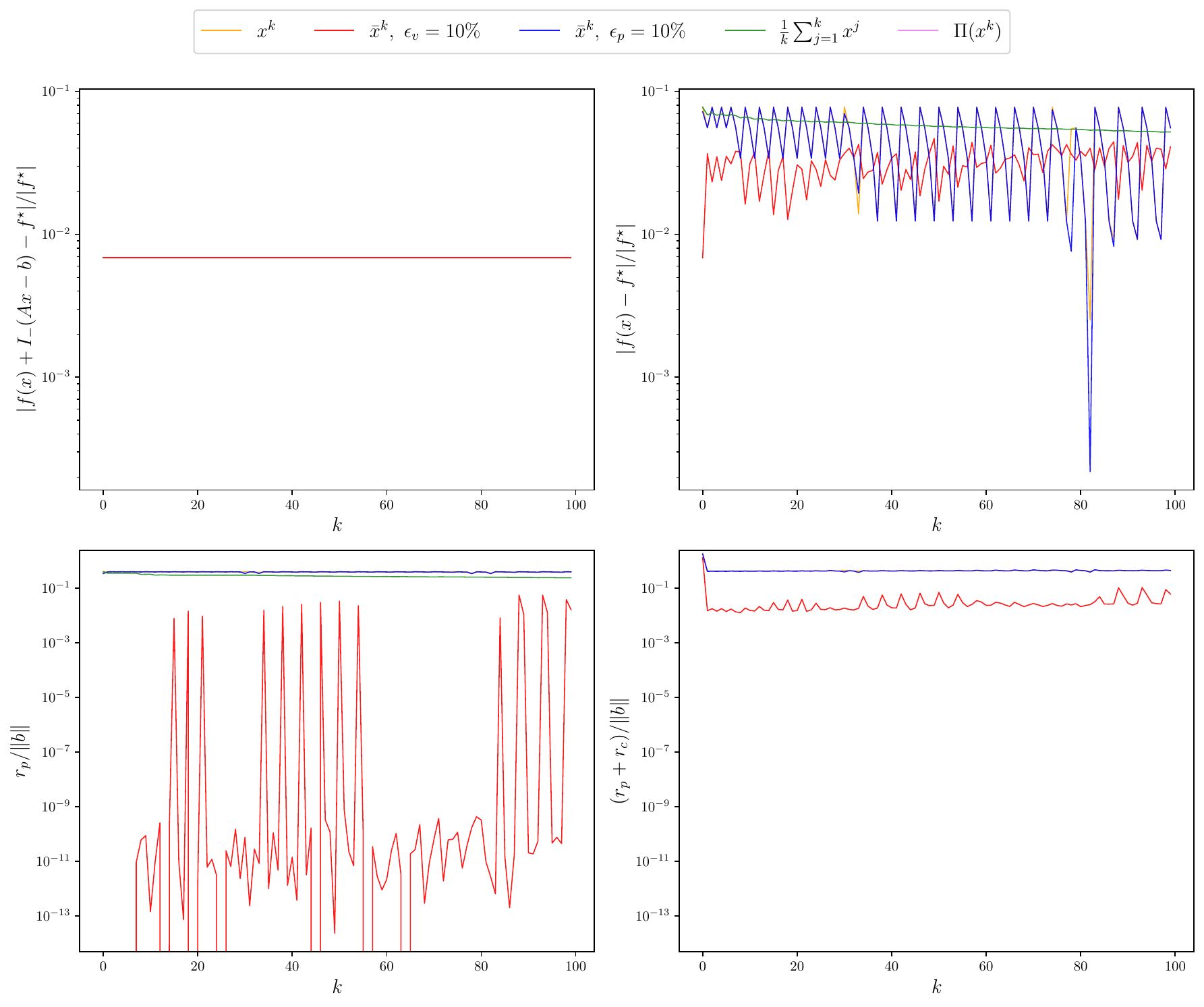}
    \caption{Best-to-date suboptimality of the most recent feasible point (top left),
    function value suboptimality (top right), relative primal violations (bottom left),
    relative residuals (bottom right)
    versus dual subgradient method iterations for the resource allocation problem.}
    \label{fig-ra-subopt-vp-subgrad}
\end{center}
\end{figure}

Similarly, Figure~\ref{fig-ra-subopt-vp-subgrad} illustrates the suboptimality and relative residuals achieved with the dual subgradient method.
MARA using a value suboptimality oracle with $\epsilon_v=10\%$ produces a primal feasible point at the first iteration
with $0.7\%$ suboptimality of the function value. However, MARA using a price perturbation oracle with $\epsilon_p=10\%$ failed to find a feasible solution within 100 iterations. 
The competing methods never recover a primal feasible point within 100 iterations, with
relative infeasibility ranging from $53\%$ to $90\%$
and
suboptimality ranging from $24\%$ to $34\%$.

\subsection{Assignment problem}\label{sec-assignment-problem}
In this example, we solve the convex relaxation of the nonconvex
problem.
Specifically, we examine the assignment problem, which involves assigning 
different projects to 
teams on the same network in a way that minimizes the total cost
while respecting the capacity constraints of
the teams.
We consider a network that is a complete bipartite graph with
each edge $(i,j)$ connecting one of the projects $i=1, \ldots, n$ to 
each of the teams $j=1, \ldots, m$.
We let individual projects and teams be agents, \ie, the
total number of agents is $K=n+m$.
Then the distributed optimization problem is given by
\[
\begin{array}{ll}
\mbox{minimize} & \sum_{i=1}^n \pi_i(x_i) + \sum_{j=1}^m \tau_j(c_j)\\ 
\mbox{subject to} & \sum_{i=1}^n (x_i)_j \leq c_j, \quad j=1, \ldots, m,
\end{array}
\]
where $\pi_i:\reals^{m}_+ \to \reals \cup \{\infty\}$ 
is the objective of the project agent $i$,
$\tau_j:\reals_+ \to \reals \cup \{\infty\}$ 
is the objective of the team agent $j$,
and $x_i \in \reals^{m}_+$, $c_j>0$ are 
variables for all $i=1, \ldots, n$ and $j=1, \ldots, m$.

We define $\pi_i(x_i)$ as the 
optimal value to the problem
\BEQ\label{e-numerical-project-obj}
\begin{array}{ll}
\mbox{minimize} & - r_i^T\tilde x_i\\ 
\mbox{subject to} & q_i \tilde x \leq x_i \\
& \tilde x_i^T \ones \leq 1 \\
& \tilde x_i \in \{0, 1\}^m,
\end{array}
\EEQ
where $\tilde x_i$ is a variable, 
$r_i \in \reals^{m}_+$ is a vector with rewards for assigning 
project $i$ to each of the $m$ teams,
and $q_i>0$ is the consumed capacity for assigning project $i$.  
The last two constraints correspond to assigning a single project to 
at most one team.

We define $\tau_j(c_j)$ as the 
optimal value to the problem
\[
\begin{array}{ll}
\mbox{minimize} & a_j(\tilde c_j - d_j)_+^2\\ 
\mbox{subject to} & c_j \leq \tilde c_j,
\quad \tilde c_j \geq 0
\end{array}
\]
where $\tilde c_j$ is a variable, 
$a_i, d_j>0$. Thus $\tau_j(c_j)$ is an increasing convex function.

\paragraph{Problem instance}
We choose the entries of $r_i$ to be uniform in $[0,1]$. 
The values $q_i$ are sampled randomly from $\{1,\ldots, \lceil n/2m \rceil\}$.  
The values $a_j$ are sampled uniformly from
$[0.8w, 1.25w]$ for $w=(\ones^TP\ones)/(30q^T\ones)$.
We sample the values $d_j$ randomly from 
$\{1, \ldots, \lceil q^T\ones/n \rceil \}$.
For this specific instance, we consider $n=200$ projects 
and $m=50$ teams.
The dimension of primal variable is $10050$
and the dual variable has dimension $50$.

To obtain convex relaxation of the assignment problem, we replace
the integer constraints in \eqref{e-numerical-project-obj}
with the box constraints $\tilde x_i \in [0, 1]^m$
for all $i=1, \ldots, m$.
This corresponds to taking the convex closure of the objective function in
assignment problem.
Since the problem is convex, we proceed with a convex version of the 
primal recovery method.

\begin{figure}
\begin{center}
\includegraphics[width=0.8\textwidth]
    {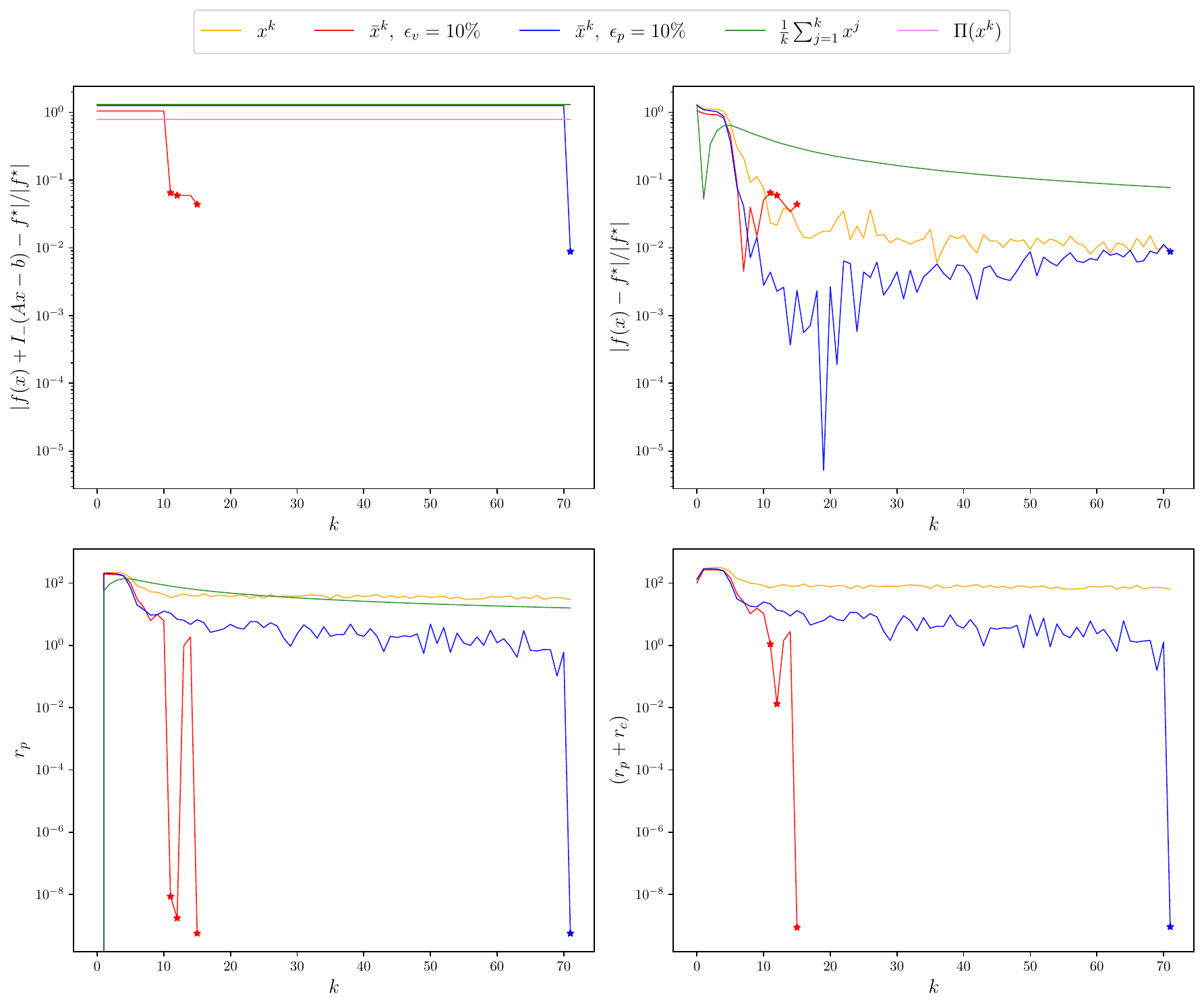}
    \caption{Best-to-date suboptimality of the most recent feasible point (top left),
    function value suboptimality (top right), primal violations (bottom left),
    residuals (bottom right)
    versus localization method iterations for
    convex relaxation of the assignment problem.}
    \label{fig-cvx-ap-subopt-vp}
\end{center}
\end{figure}

\paragraph{Results}
Figure~\ref{fig-cvx-ap-subopt-vp} illustrates the suboptimality and
relative residuals achieved with the localization method
for the assignment problem.
At iteration $15$, the $\epsilon_v=10\%$ primal recovery method  finds a 
feasible primal
point $\bar x^k$ with $4.4\%$ suboptimality of the function value.
At iteration $71$, the $\epsilon_p=10\%$ primal recovery method finds a 
feasible primal
point $\bar x^k$ with a function value suboptimality of $0.9\%$.
In contrast, the primal points $x^k$, primal average points $\frac{1}{k}\sum_{j=1}^k x^j$,
and the projected points $\Pi(x^k)$ are
primal feasible only at the initial iteration, with
suboptimalities of
$130\%$, $130\%$ and $79\%$, respectively.

\begin{figure}
\begin{center}
\includegraphics[width=0.8\textwidth]
    {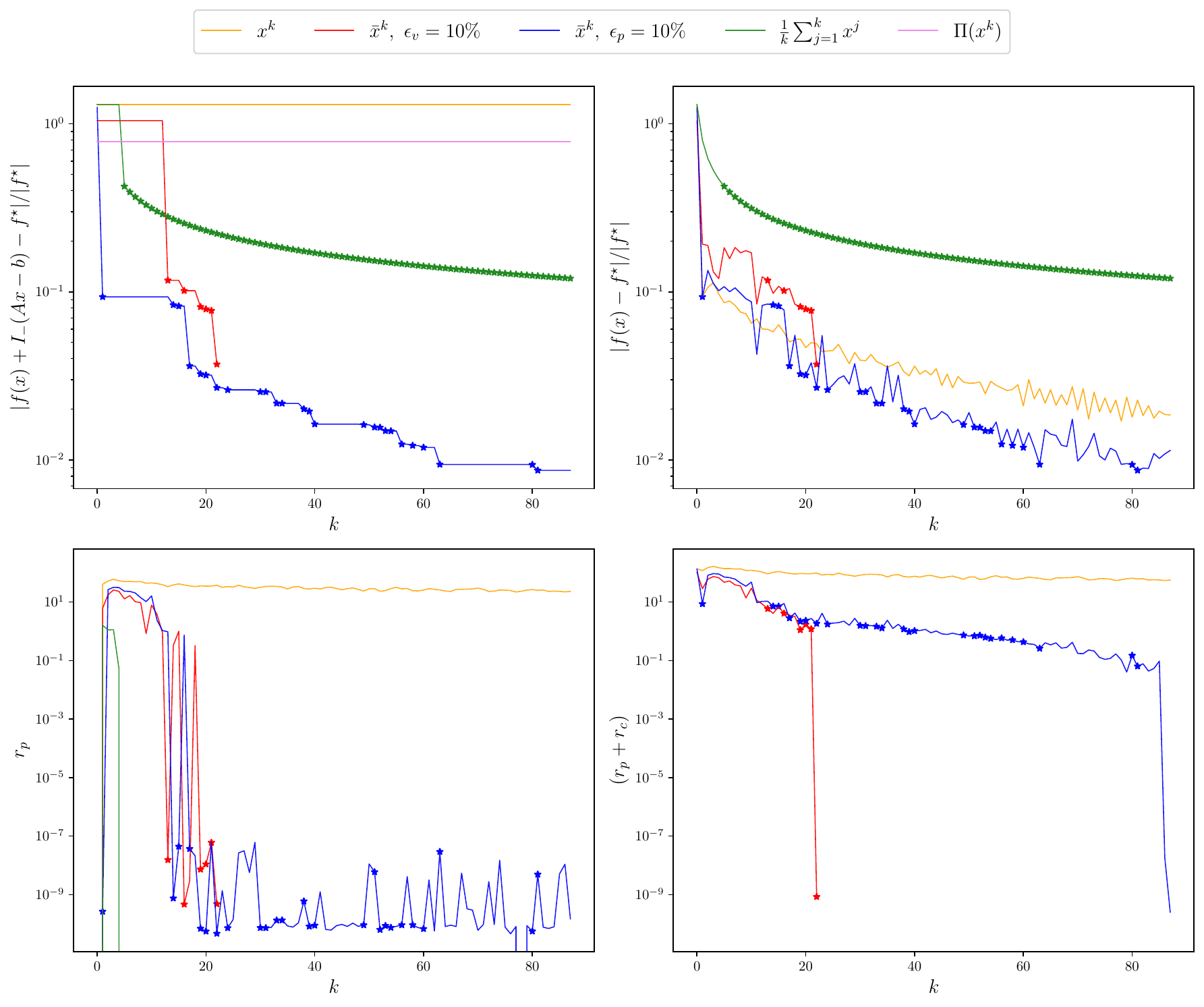}
    \caption{Best-to-date suboptimality of the most recent feasible point (top left),
    function value suboptimality (top right), primal violations (bottom left),
     residuals (bottom right)
    versus dual subgradient method iterations for
    convex relaxation of the assignment problem.}
    \label{fig-cvx-ap-subopt-vp-subgrad}
\end{center}
\end{figure}

Similarly, Figure~\ref{fig-cvx-ap-subopt-vp-subgrad} illustrates suboptimality and
relative residuals achieved with the dual subgradient method.
Although all of the methods are primal feasible at initialization,
MARA  achieve better suboptimality quickly as well.
In particular, MARA with $\epsilon_v=10\%$ and $\epsilon_p=10\%$ 
produces primal feasible points at iterations $22$ and $20$
with  function value suboptimality of $3.7\%$ and $3.2\%$, respectively.
In contrast, the competing methods have
suboptimality ranging from $23\%$ to $130\%$.

\subsection{Multi-commodity flow problem}\label{sec-mcf-problem}
We consider a variation of the multi-commodity flow problem from 
\cite[\S4.3]{parshakova2024implementation}.
Let $A\in \reals^{p \times q}$ be the incidence matrix
of a graph with $p$ vertices or nodes
and $q$ directed edges.
The network supports the flow of $K$ different commodities.
Each commodity $i=1, \ldots, K$ is shipped from a source node 
$r_i \in \{1,\ldots, p\}$
to a sink $s_i \in \{1,\ldots, p\}$ over the network edges
with total flow given by $d_i\geq 0$.
We also define $z_i\in \reals_+^q$ as the vector of flows of commodity $i$ 
on the edges.
Then the flow conservation is given by
\[
Az_i +  d_i(e_{r_i} - e_{s_i})= 0, \quad i=1, \ldots, K,
\]
where 
$e_j$ is the $j$th unit vector.

We consider the multi-commodity problem 
\[
\begin{array}{ll} \mbox{maximize} & \sum_{i=1}^K b_i \sqrt{d_i}\\
\mbox{subject to} & 0 \leq z_i \leq x_i, \quad i=1, \ldots, K\\
& Az_i+d_i(e_{r_i} - e_{s_i}) =0, \quad i=1, \ldots, K\\
& v_1 x_1 + \cdots + v_K x_K \leq c, \quad x_i \geq 0, \quad i=1, \ldots, K,
\end{array}
\]
where $z_i$, $d_i$, and $x_i$ are variables, 
$c \in \reals_+^q$ is a given vector of edge capacities,
and $v\in \reals_+^K$ is a given vector of commodity volumes.

We define the agent function 
$f_i(x_i)$
as the optimal value of the single commodity problem
\[
\begin{array}{ll} \mbox{minimize} & -b_i \sqrt{d_i}\\
\mbox{subject to} & 0 \leq z_i \leq x_i\\
& Az_i+d_i(e_{r_i} - e_{s_i}) =0 \\
& x_i \leq R \ones, 
\end{array}
\]
with variables $z_i$ and $d_i$, 
for all $i=1, \ldots, K$.
Then the distributed optimization problem takes the form
\[
\begin{array}{ll}
\mbox{minimize} & \sum_{i=1}^K f_i(x_i)\\ 
\mbox{subject to} & v_1 x_1 + \cdots + v_K x_K \leq c.
\end{array}
\]

\paragraph{Problem instance}
We consider an example with $K=100$ commodities,
and a graph with $p=15$ nodes and $q=100$ edges. Edges are generated
randomly from pairs of nodes, with an additional cycle passing through 
all vertices (to ensure that the graph is strongly connected,
\ie, there is a directed path from any node to any other). We randomly
choose the source-destination pairs $(r_i, s_i)$.
We choose commodity capacities $\tilde c_j$ from a uniform distribution on $[0.2,2]$,
and $b_i$ are chosen uniformly on $[0.5,1.5]$.
Volumes are sampled from $\exp N(0,1)$. Then the capacities are set as 
expected volume per edge, \ie, 
$c_j = \tilde c_j \exp(1/2)$ for all $j=1, \ldots, m$. 
The dimension of primal variable is $10000$
and of dual variable is $100$.

\begin{figure}
\begin{center}
    \includegraphics[width=0.8\textwidth]{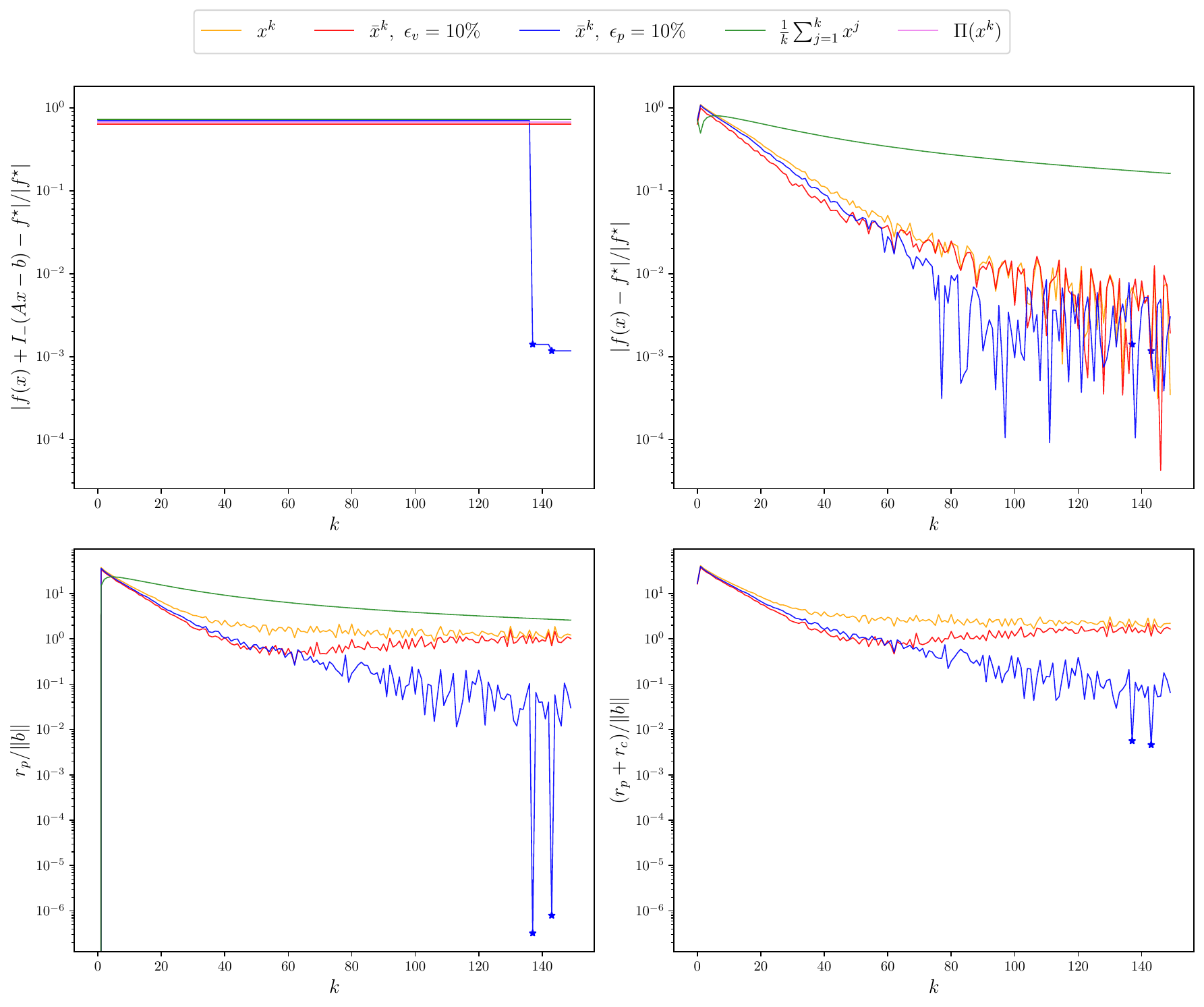}
    \caption{Best-to-date suboptimality of the most recent feasible point (top left),
    function value suboptimality (top right), relative primal violations (bottom left),
    relative residuals (bottom right)
    versus localization method iterations for the
    multi-commodity problem. }
    \label{fig-mcf-subopt-vp}
\end{center}
\end{figure}

\paragraph{Results}
Figure~\ref{fig-mcf-subopt-vp} illustrates suboptimality and
relative residuals achieved with the localization method
for the multi-commodity flow problem.
At iteration $62$, MARA with 
$\epsilon_v=10\%$
finds a primal
point $\bar x^k$ with $28\%$ relative primal infeasibility, 
$1.8\%$ function value suboptimality,
and $2.8\times$ smaller primal residual than $x^k$.
At iteration $137$ MARA with 
$\epsilon_p=10\%$
finds a primal
feasible point with 
$0.1\%$ function value suboptimality.

In contrast, at iteration $143$, $x^k$ is primal infeasible 
with $79\%$ relative infeasibility
and $0.1\%$ suboptimality.
At iteration $149$, the primal average value is infeasible with
$258\%$ relative infeasibility
and $16\%$ suboptimality.
Moreover, except for the first iteration, the projected point $\Pi(x^k)$
is everywhere outside the domain of the agents'
functions.
The competing methods are primal feasible only at the initial iteration with the following suboptimalities:
the primal point $x^k$ is $72\%$ suboptimal, 
primal average point $\frac{1}{k}\sum_{j=1}^k x^j$ is $72\%$ suboptimal,
and the projected points $\Pi(x^k)$ is $67\%$ suboptimal.

\begin{figure}
\begin{center}
    \includegraphics[width=0.8\textwidth]{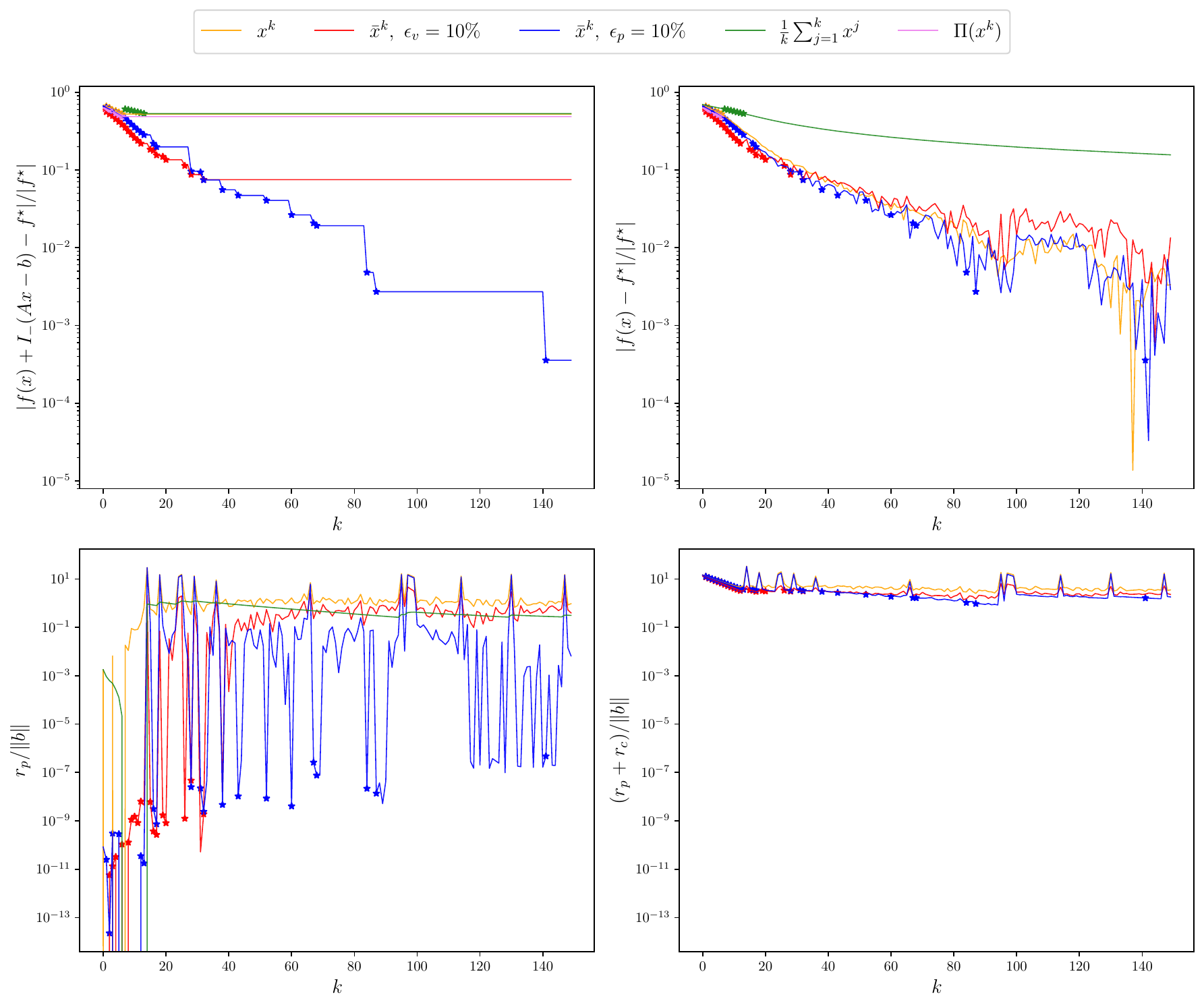}
    \caption{Best-to-date suboptimality of the most recent feasible point (top left),
    function value suboptimality (top right), relative primal violations (bottom left),
    relative residuals (bottom right)
    versus dual subgradient method iterations for the
    multi-commodity problem. }
    \label{fig-mcf-subopt-vp-subgrad}
\end{center}
\end{figure}

Similarly, Figure~\ref{fig-mcf-subopt-vp-subgrad} illustrates suboptimality and 
relative residuals achieved with the dual subgradient method.
MARA with both approximate oracles recovers primal feasible points with suboptimality
levels that are a factor
$\times 5$ to $\times 7$ lower than those of the competing methods.

\subsection{Shipment problem}\label{sec-shipment-problem}
We consider a shipment problem given by the following capacitated optimal transport formulation
\[
\begin{array}{ll}
\mbox{minimize} & \Tr(C^T X) \\ 
\mbox{subject to} & X \ones = \mu_s \\
& X_{i,j} \geq 0, \quad i = 1, \ldots, K, \quad j = 1, \ldots, m \\
& X^T \ones = \mu_t \\
& X^T v \leq c
\end{array}
\]
where $C:\reals^{K \times m}_+$ is a cost matrix,
$\mu_s\in \reals^K$ is a source distribution, 
$\mu_t \in \reals^m$ is a target distribution, 
$v \in \reals^K$ is vector with volumes, 
$c \in \reals^m$ is vector with capacities, 
and $X \in \reals^{K \times m}_+$ is 
a variable. 

Let $c_i$ and $x_i$ be $i$th row of $C$ and $X$ respectively.
We define the agent function $f_i:\reals^m \to \reals\cup \{\infty \}$ as
\[
f_i(x_i)=  
\begin{cases}
    c_i^T x_i & x_i \geq 0, \quad x_i^T \ones = (\mu_s)_i \\
    \infty & \text{otherwise}
\end{cases}
\]
for all $i=1, \ldots, K$.
Then the distributed optimization problem takes the form
\[
\begin{array}{ll}
\mbox{minimize} & \sum_{i=1}^K f_i(x_i)\\ 
\mbox{subject to} & X^T \ones = \mu_t \\
& X^T v \leq c \\
& X = (x_1, \ldots, x_K).
\end{array}
\]

\paragraph{Problem instance}
We sample entries of source and target locations, $s\in \reals^{K \times d}$ and
$t\in \reals^{m \times d}$,
from the standard normal.
Then entries of cost matrix are defined as $C_{ij}=\|s_i - t_j\|_2$
for all $i = 1, \ldots, K$ and $j = 1, \ldots, m$.
The source and target distributions $\mu_s$ and $\mu_t$ respectively are sampled
from log-normal and normalized to sum to one.
Volumes are sampled from $\exp N(0,\sigma^2)$. Then capacities are set as 
expected volume per target location, \ie, 
$c_j = b_j \exp(\sigma^2/2)$ for all $j=1, \ldots, m$. 

For this specific instance, we consider $m=25$,
$K=100$ agents, $d=10$ and $\sigma=0.8$.
The dimension of primal variable is $2500$
and the dual variable has dimension $50$.

\begin{figure}
\begin{center}
    \includegraphics[width=0.8\textwidth]{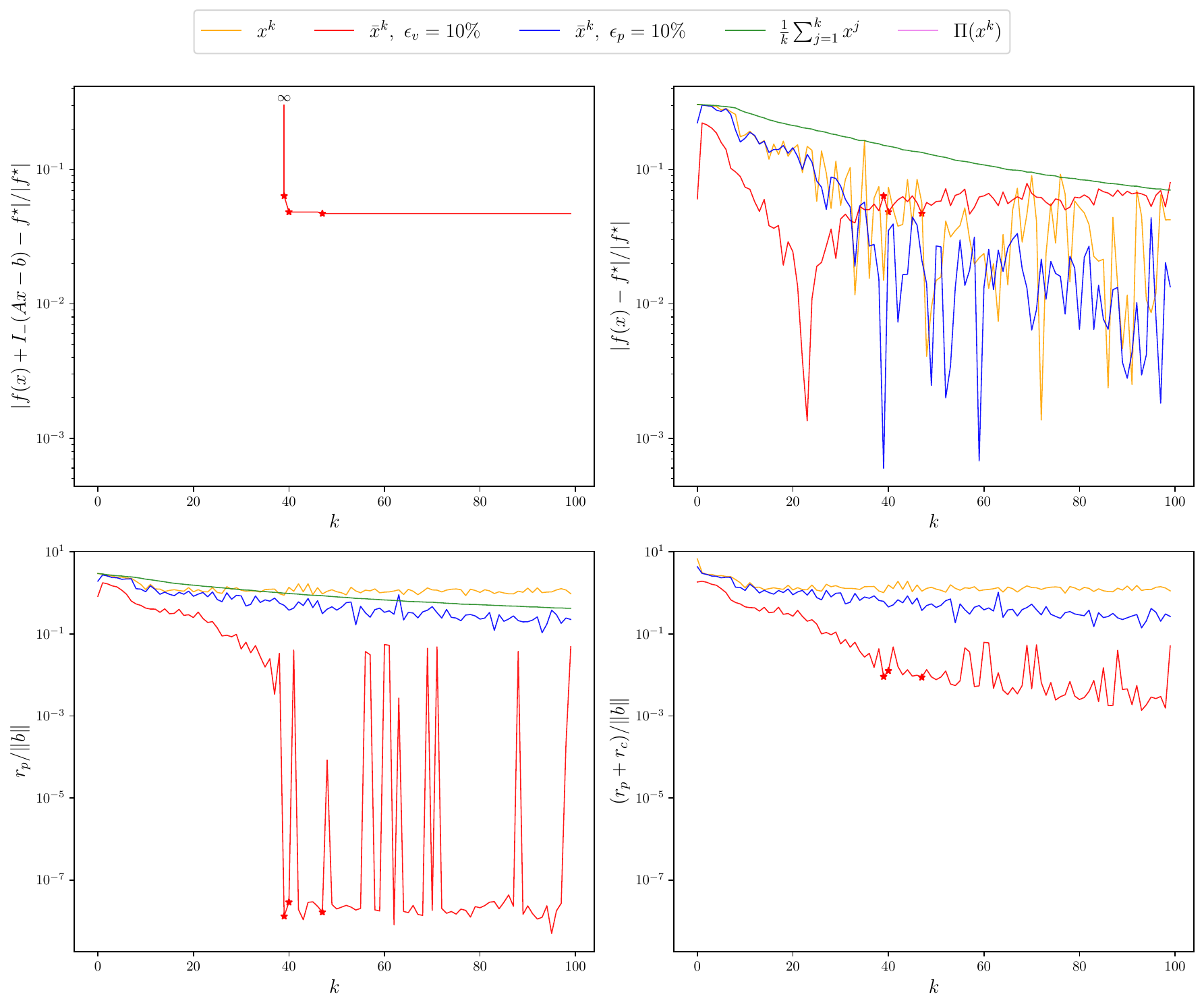}
    \caption{Best-to-date suboptimality of the most recent feasible point (top left),
    function value suboptimality (top right), relative primal violations (bottom left),
    relative residuals (bottom right)
    versus localization method iterations for the
    shipment problem. }
    \label{fig-ot-subopt-vp}
\end{center}
\end{figure}

\paragraph{Results}
Figure~\ref{fig-ot-subopt-vp} illustrates suboptimality and relative residuals 
achieved with localization method
for shipment problem.
At iteration $39$, MARA with 
$\epsilon_v=10\%$
finds a feasible primal
point $\bar x^k$ with $6.3\%$ function value suboptimality.
At iteration $93$ MARA with 
$\epsilon_p=10\%$
finds a primal
point $\bar x^k$ with $0.3\%$ function value suboptimality
and $8\times$ reduction of primal violation relative to $x^k$
and $4\times$ reduction relative to the primal average point.

In contrast, at iteration $77$,
$x^k$ is primal infeasible with $86\%$ relative infeasibility
and $6.6\%$ suboptimality.
At iteration $99$, the primal average point is infeasible 
with $42\%$ relative infeasibility
and $7\%$ suboptimality.
The projected point $\Pi(x^k)$ is everywhere 
outside the domain of the agents' objective
functions.

\begin{figure}
\begin{center}
    \includegraphics[width=0.8\textwidth]{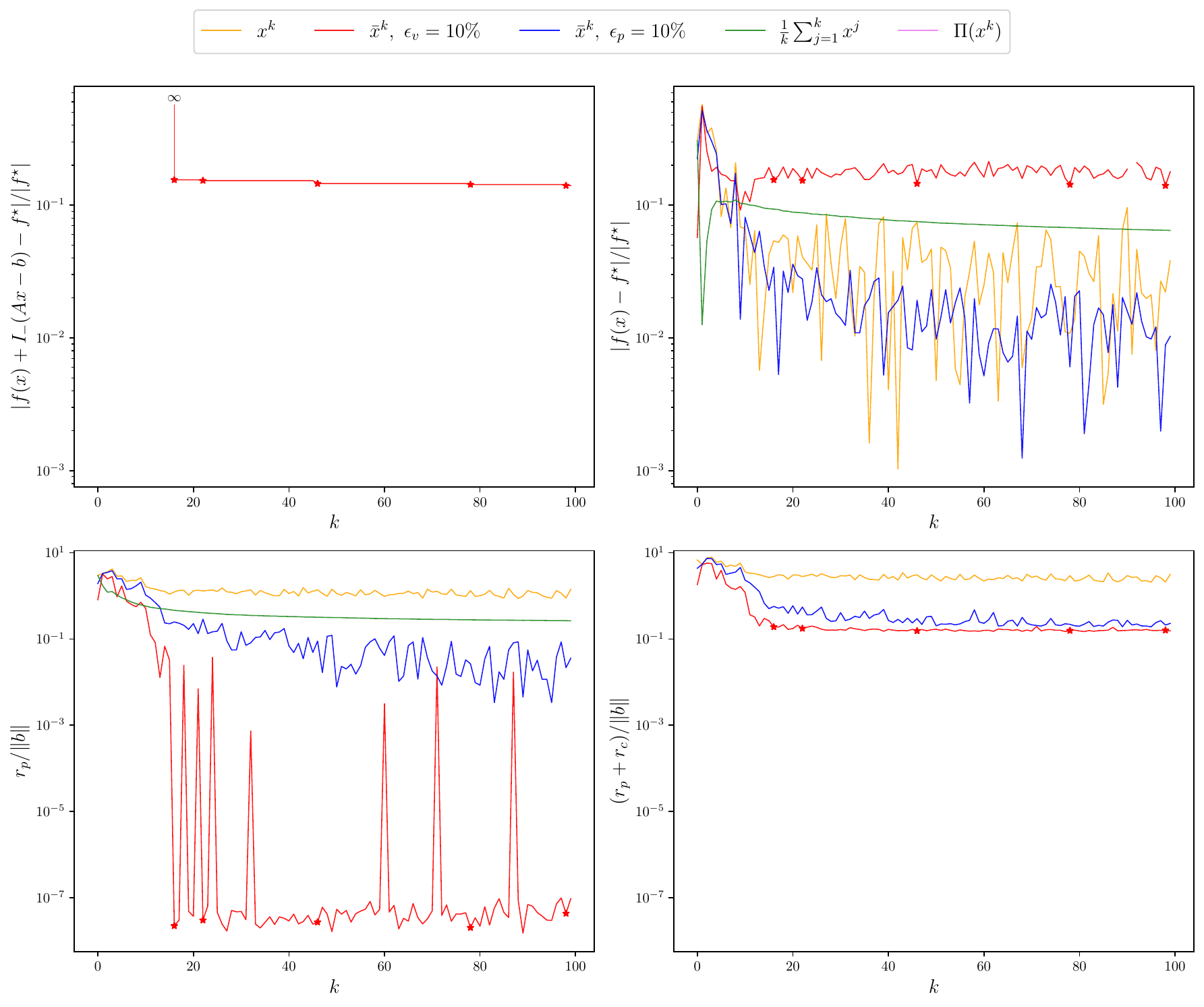}
    \caption{Best-to-date suboptimality of the most recent feasible point (top left),
    function value suboptimality (top right), relative primal violations (bottom left),
    relative residuals (bottom right)
    versus dual subgradient method iterations for the
    shipment problem. }
    \label{fig-ot-subopt-vp-subgrad}
\end{center}
\end{figure}

Similarly, Figure~\ref{fig-ot-subopt-vp-subgrad} illustrates suboptimality and
relative residuals achieved with the dual subgradient method.
MARA with $\epsilon_v=10\%$ recovers a primal feasible point
at iteration $16$ with suboptimality of $15.5\%$, while the MARA price method with $\epsilon_p=10\%$ fails to find a feasible solution.
All competing methods are unable to recover a primal feasible point, with
relative infeasibility ranging from $26\%$ to $87\%$
and
suboptimality ranging from $1.1\%$ to $6.5\%$.

\subsection{Approximation trade-offs}\label{sec-exp-tradeoffs}

Next we evaluate the performance of MARA using the price localization method
by varying two factors:
1) oracle suboptimality, $\epsilon=1\%$ and $\epsilon=10\%$,
and 2) number of responses, $N=10$ and $N=50$.

\paragraph{Oracle suboptimality $\epsilon$}
Recall from \S\ref{sec-approx-oracle}, that the oracle's suboptimality has 
additive error term proportional to $\epsilon$.
Thus, we expect smaller values of $\epsilon$
to improve the function value suboptimality of MARA.
At the same time primal points generated by the approximate oracle with smaller 
$\epsilon$
are a subset of those with larger $\epsilon$, which implies
slower reduction in residuals.

\begin{figure}[h!]
\begin{center}
   \includegraphics[width=0.8\textwidth]
    {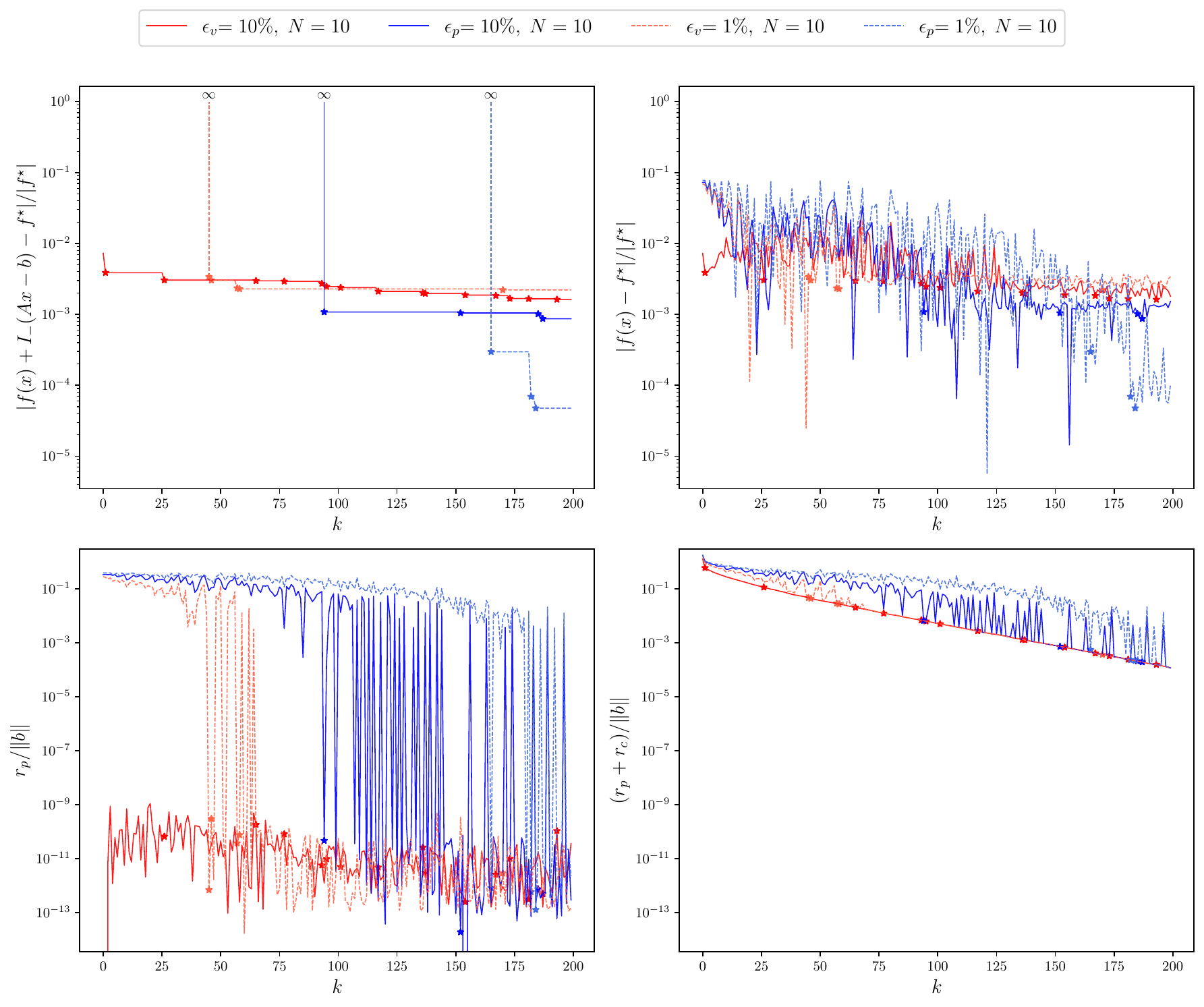}
    \caption{Best-to-date suboptimality of the most recent feasible point (top left),
    function value suboptimality (top right), relative primal violations (bottom left),
    relative residuals (bottom right)
    versus localization method iterations  for
    the resource allocation problem with  $\epsilon=\{1\%, 10\%\}$
    and $N=10$. }
    \label{fig-ra-eps-vary}
\end{center}
\end{figure}

\begin{figure}
\begin{center}
    \includegraphics[width=0.8\textwidth]
    {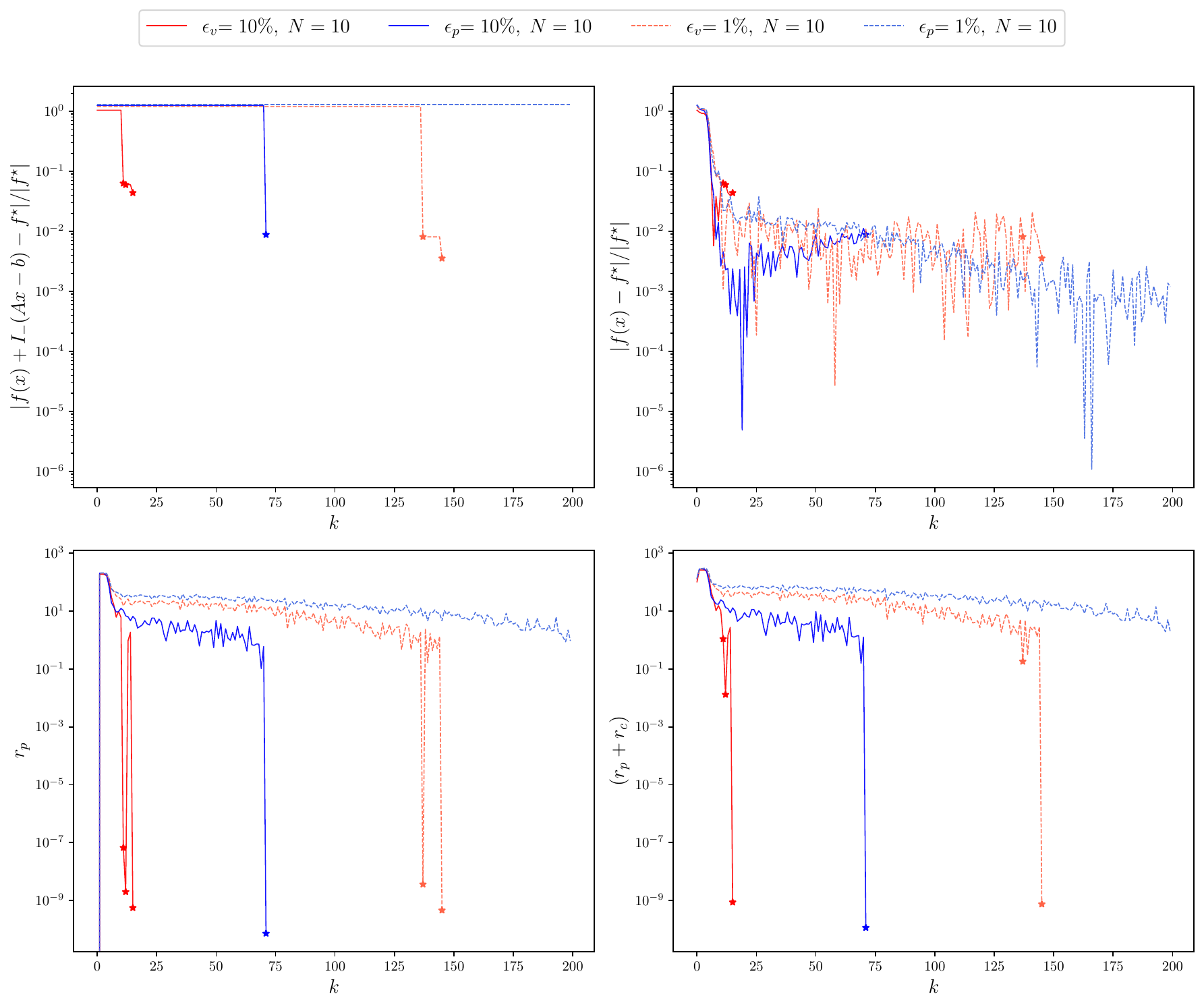}
    \caption{Best-to-date suboptimality of the most recent feasible point (top left),
    function value suboptimality (top right), primal violations (bottom left),
    residuals (bottom right)
    versus localization method iterations for
    convex relaxation of the assignment problem for  $\epsilon=\{1\%, 10\%\}$ and $N=10$. }
    \label{fig-ap-eps-vary}
\end{center}
\end{figure}

\begin{figure}[h!]
\begin{center}
   \includegraphics[width=0.8\textwidth]
    {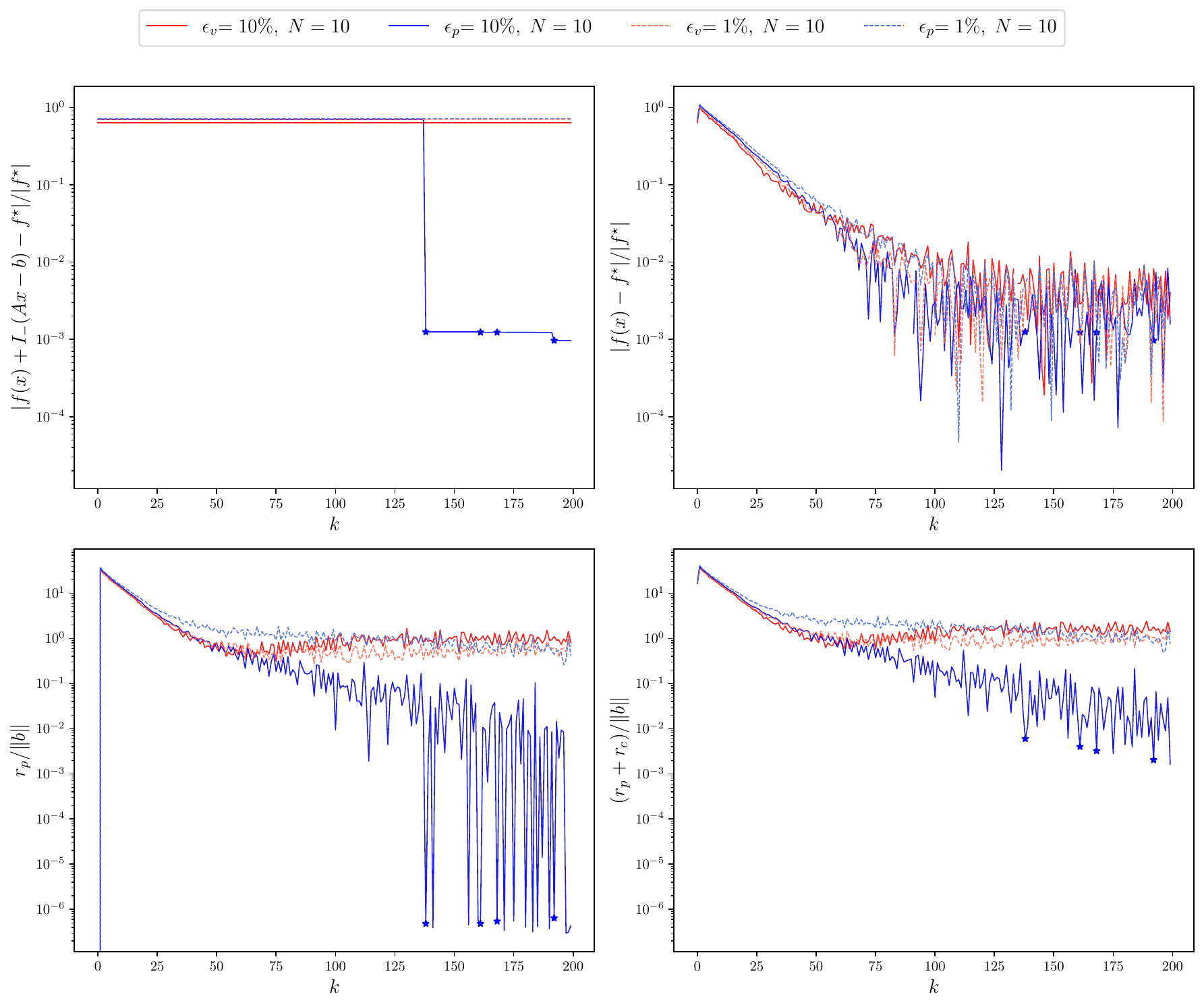}
    \caption{Best-to-date suboptimality of the most recent feasible point (top left),
    function value suboptimality (top right), relative primal violations (bottom left),
    relative residuals (bottom right)
    versus localization method iterations  for
    the multi-commodity flow problem with  $\epsilon=\{1\%, 10\%\}$ and $N=10$. }
    \label{fig-mcf-eps-vary}
\end{center}
\end{figure}

\begin{figure}[h!]
\begin{center}
   \includegraphics[width=0.8\textwidth]
    {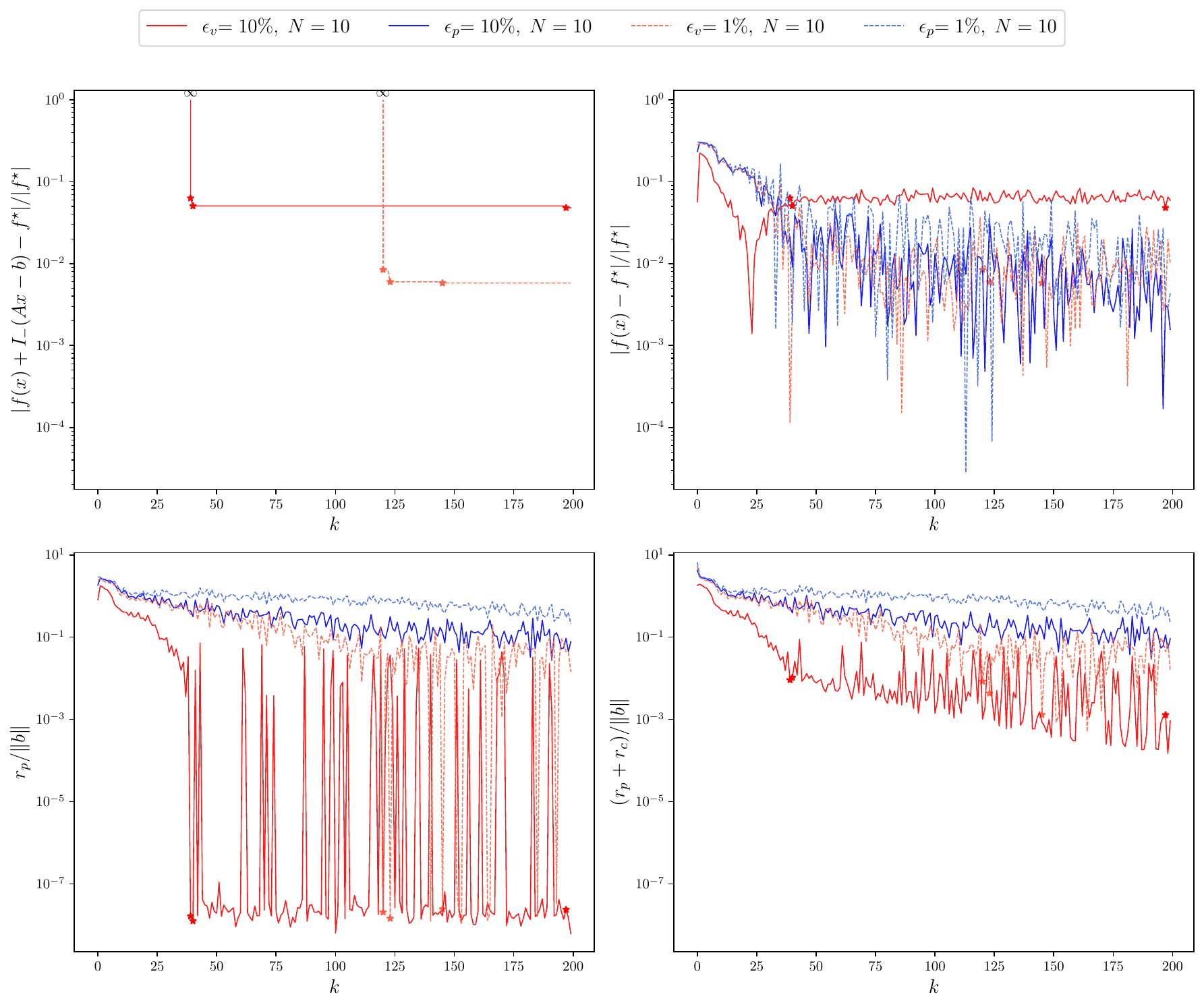}
    \caption{Best-to-date suboptimality of the most recent feasible point (top left),
    function value suboptimality (top right), relative primal violations (bottom left),
    relative residuals (bottom right)
    versus localization method iterations  for
    the shipment problem with  $\epsilon=\{1\%, 10\%\}$ and $N=10$. }
    \label{fig-ot-eps-vary}
\end{center}
\end{figure}

This behavior is indeed observed in our numerical experiments; see Figures \ref{fig-ra-eps-vary},
\ref{fig-ap-eps-vary}, \ref{fig-ot-eps-vary}, and \ref{fig-mcf-eps-vary}.
For example, as shown in Figure \ref{fig-ra-eps-vary}, MARA with 
$\epsilon_p=10\%$ finds a primal feasible primal point at iteration $94$, 
much faster than 
the method with $\epsilon_p=1\%$, which takes $165$ iterations.
Their best suboptimalities for primal feasible points
are $0.1\%$ and $0.005\%$ respectively.
Similarly, MARA with 
$\epsilon_v=10\%$ finds a feasible primal point at iteration $0$, 
while the method $\epsilon_v=1\%$ takes $45$ iterations.
However, the best suboptimalities of the feasible points for the two methods are 
comparable, $0.16\%$ and $0.22\%$, respectively.

\paragraph{Number of responses $N$}
As discussed in \S\ref{sec-approx-oracle}, 
the primal point returned by the primal recovery
method is a convex combination of $N_i$ $\epsilon$-suboptimal primal points
for each agent $i=1, \ldots, K$.
Increasing $N_i$ therefore enlarges the convex hull of the candidate
primal points,
implying a
faster reduction in residuals.

\begin{figure}
\begin{center}
   \includegraphics[width=0.8\textwidth]
    {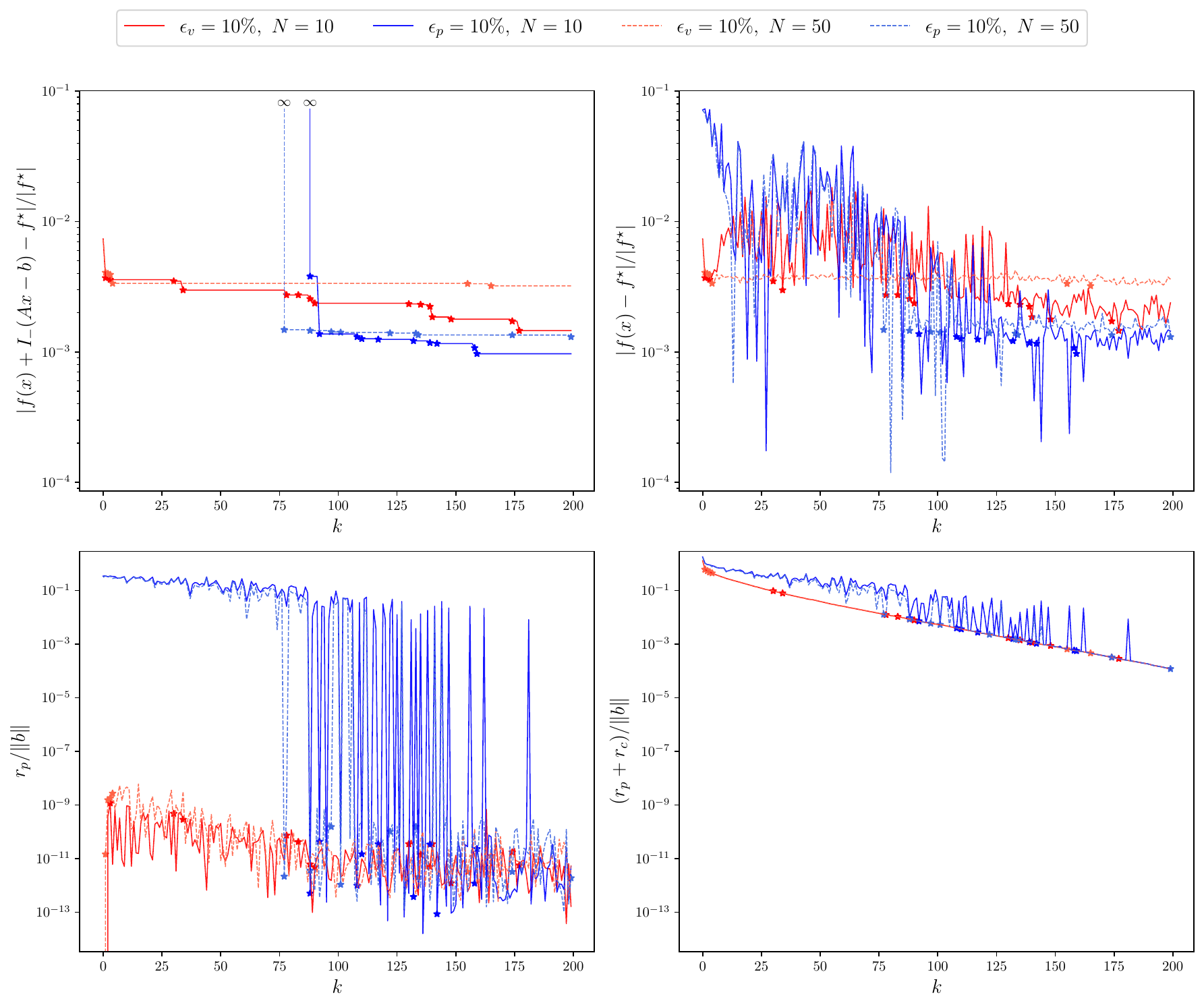}
    \caption{Best-to-date suboptimality of the most recent feasible point (top left),
    function value suboptimality (top right), relative primal violations (bottom left),
    relative residuals (bottom right)
    versus localization method iterations  for
    the resource allocation problem with $N=\{10, 50\}$ 
    and $\epsilon=10\%$. }
    \label{fig-ra-N-vary}
\end{center}
\end{figure}

\begin{figure}
\begin{center}
   \includegraphics[width=0.8\textwidth]
    {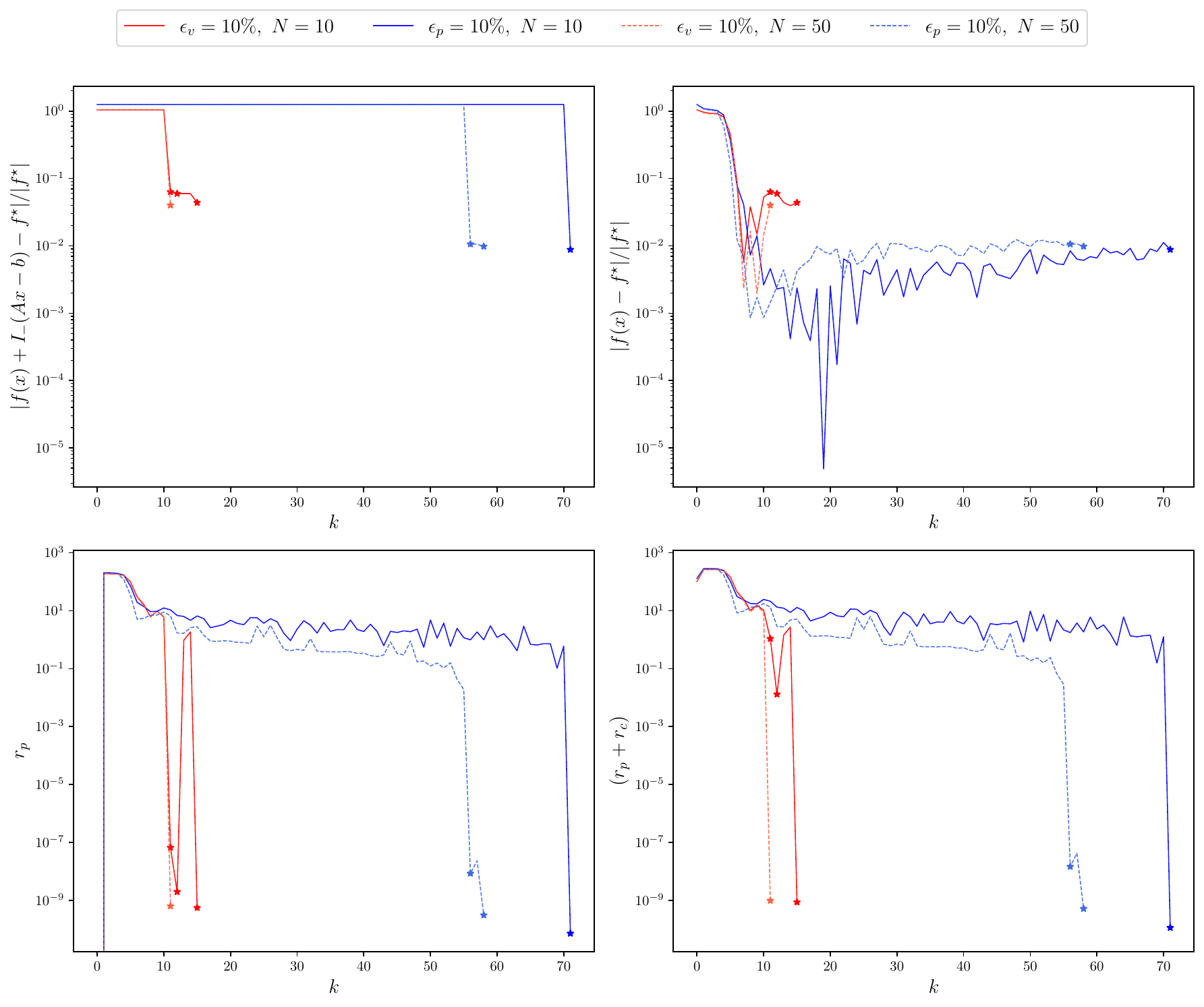}
    \caption{Best-to-date suboptimality of the most recent feasible point (top left),
    function value suboptimality (top right), relative primal violations (bottom left),
    relative residuals (bottom right)
    versus localization method iterations for
    convex relaxation of the assignment problem with $N=\{10, 50\}$ 
    and $\epsilon=10\%$. }
    \label{fig-ap-N-vary}
\end{center}
\end{figure}

\begin{figure}
\begin{center}
   \includegraphics[width=0.8\textwidth]
    {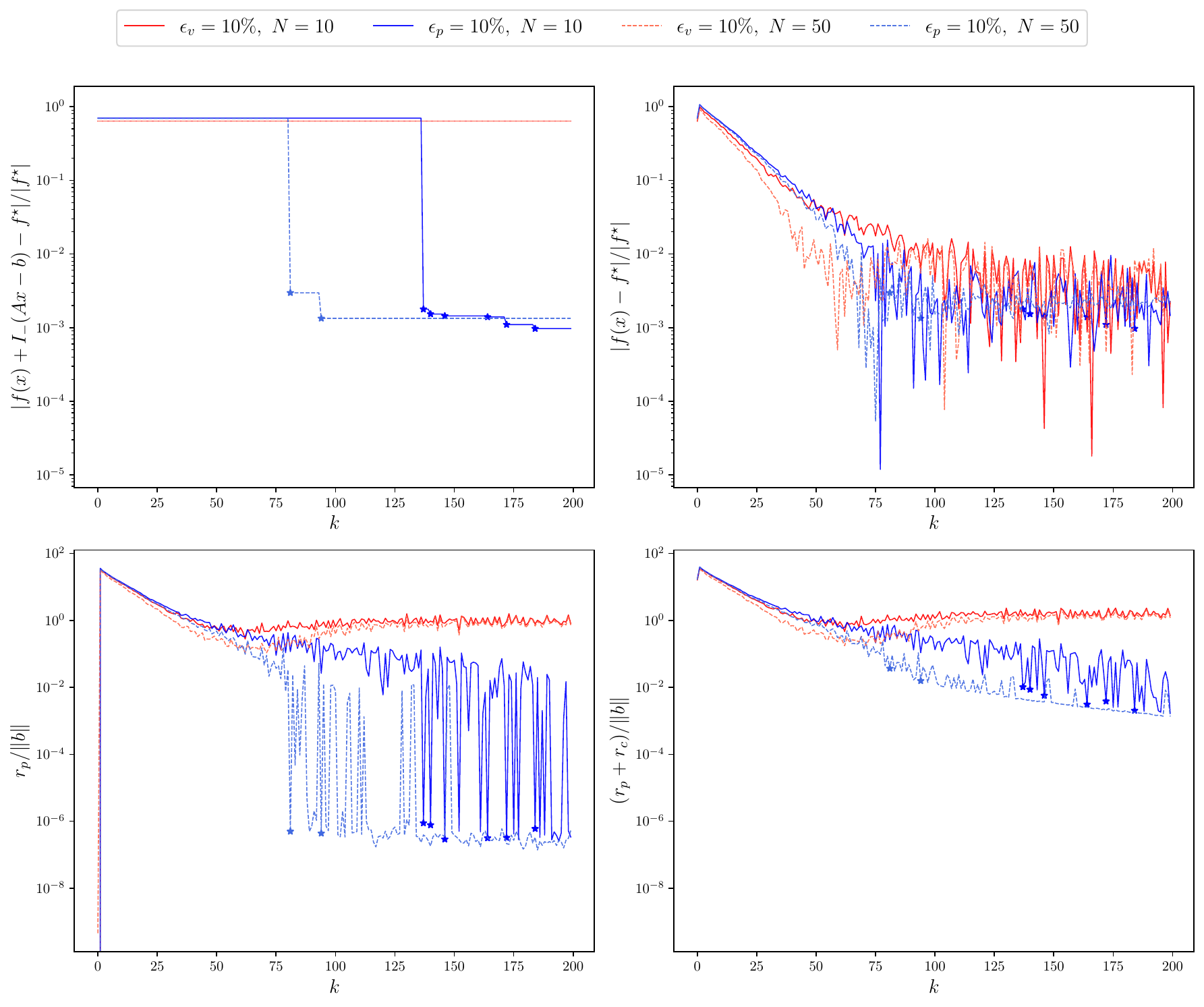}
    \caption{Best-to-date suboptimality of the most recent feasible point (top left),
    function value suboptimality (top right), relative primal violations (bottom left),
    relative residuals (bottom right)
    versus localization method iterations  for
    the multi-commodity flow problem with $N=\{10, 50\}$ 
    and $\epsilon=10\%$. }
    \label{fig-mcf-N-vary}
\end{center}
\end{figure}

\begin{figure}
\begin{center}
    \includegraphics[width=0.8\textwidth]
    {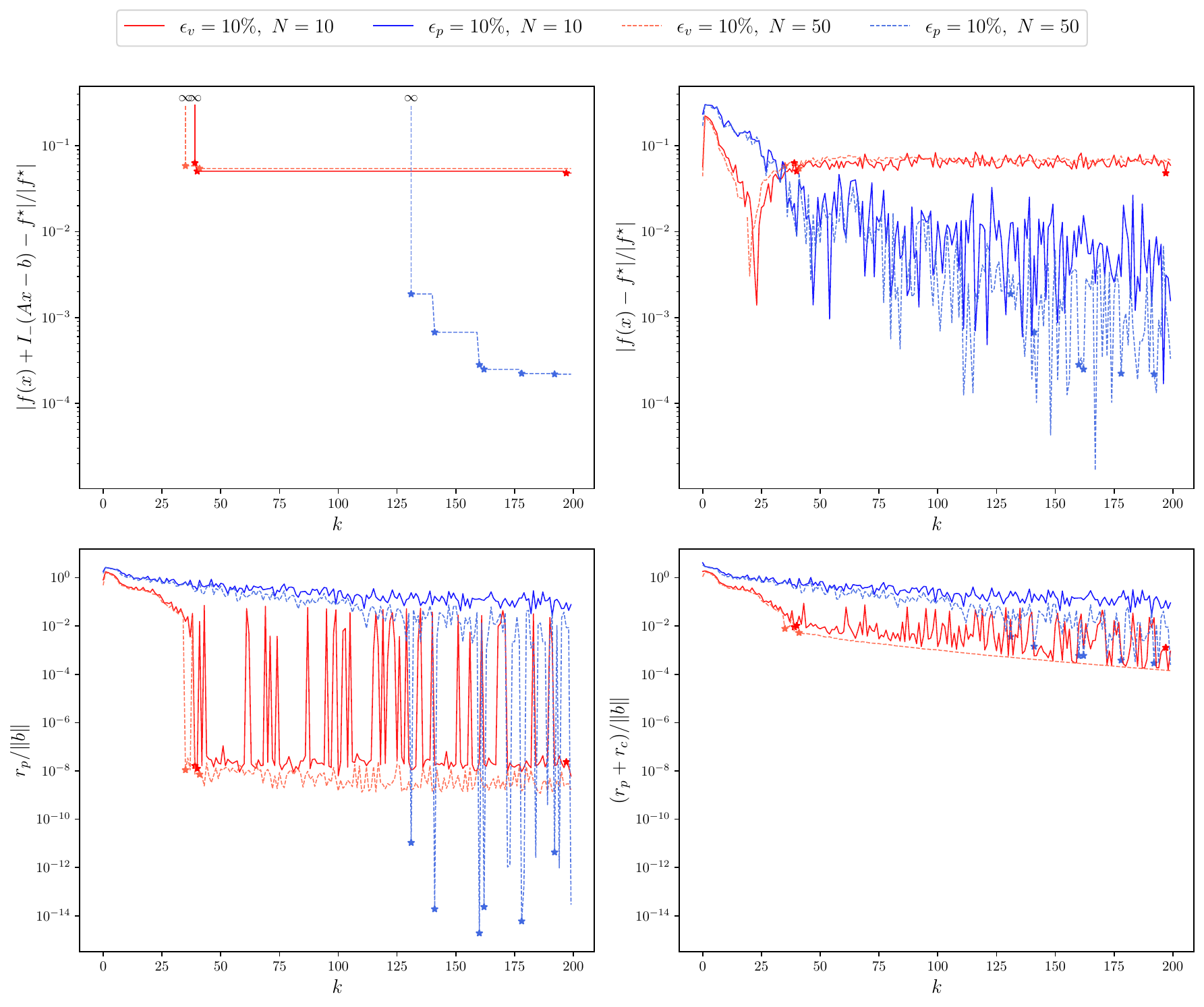}
    \caption{Best-to-date suboptimality of the most recent feasible point (top left),
    function value suboptimality (top right), relative primal violations (bottom left),
    relative residuals (bottom right)
    versus localization method iterations  for
    the shipment problem for  $N=\{10, 50\}$ 
    and $\epsilon=10\%$. }
    \label{fig-ot-N-vary}
\end{center}
\end{figure}

This behavior is observed in our numerical experiments,
see
Figures \ref{fig-ra-N-vary}, \ref{fig-ap-N-vary},
\ref{fig-ot-N-vary} and \ref{fig-mcf-N-vary}.
For example, as shown in Figure \ref{fig-mcf-N-vary}, 
(excluding the initial iteration),
MARA with 
$\epsilon_p=10\%$ and $N=10$ finds a primal feasible primal point at iteration $137$, 
much slower than
the method with $\epsilon_p=10\%$ and $N=50$, which takes $81$ iterations.
Their best suboptimalities for feasible points for two methods are 
comparable, $0.10\%$ and $0.13\%$, respectively.
The results of this section suggest it might be beneficial to split the budget of 
samples $N_i$ across oracles with different suboptimality. 

\subsection{Including history in MARA }\label{sec-exper-history}
\begin{figure}[t]
\begin{center}
    \centering
        \includegraphics[width=0.9\textwidth]{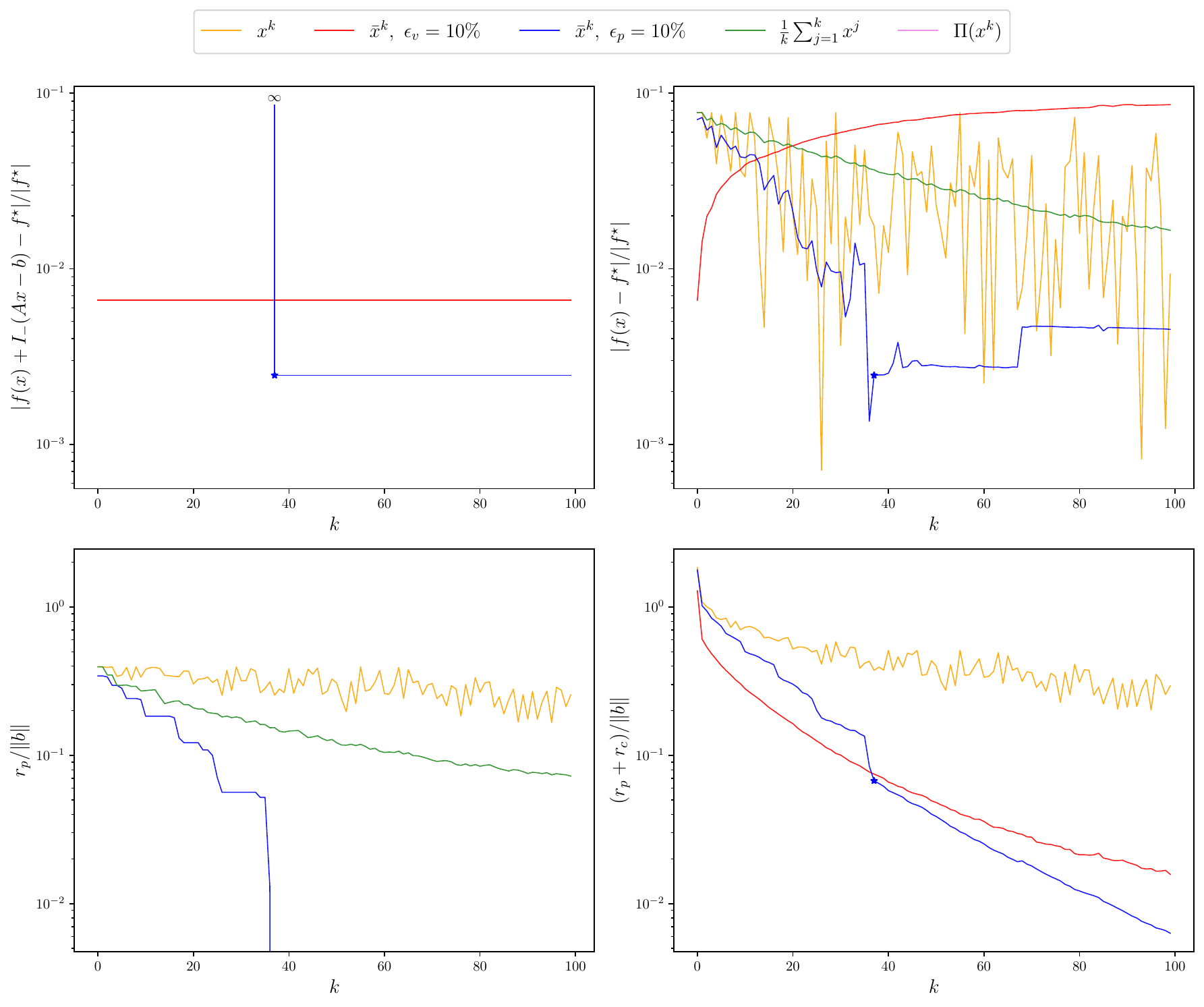}
    \caption{Best-to-date suboptimality of the most recent feasible point (top left),
    function value suboptimality (top right), relative primal violations (bottom left),
    relative residuals (bottom right)
    versus localization method iterations for the resource allocation problem
    with history $H=3$.}
    \label{fig-ra-subopt-vp-h}
\end{center}
\end{figure}

\begin{figure}
\begin{center}
\includegraphics[width=0.8\textwidth]
    {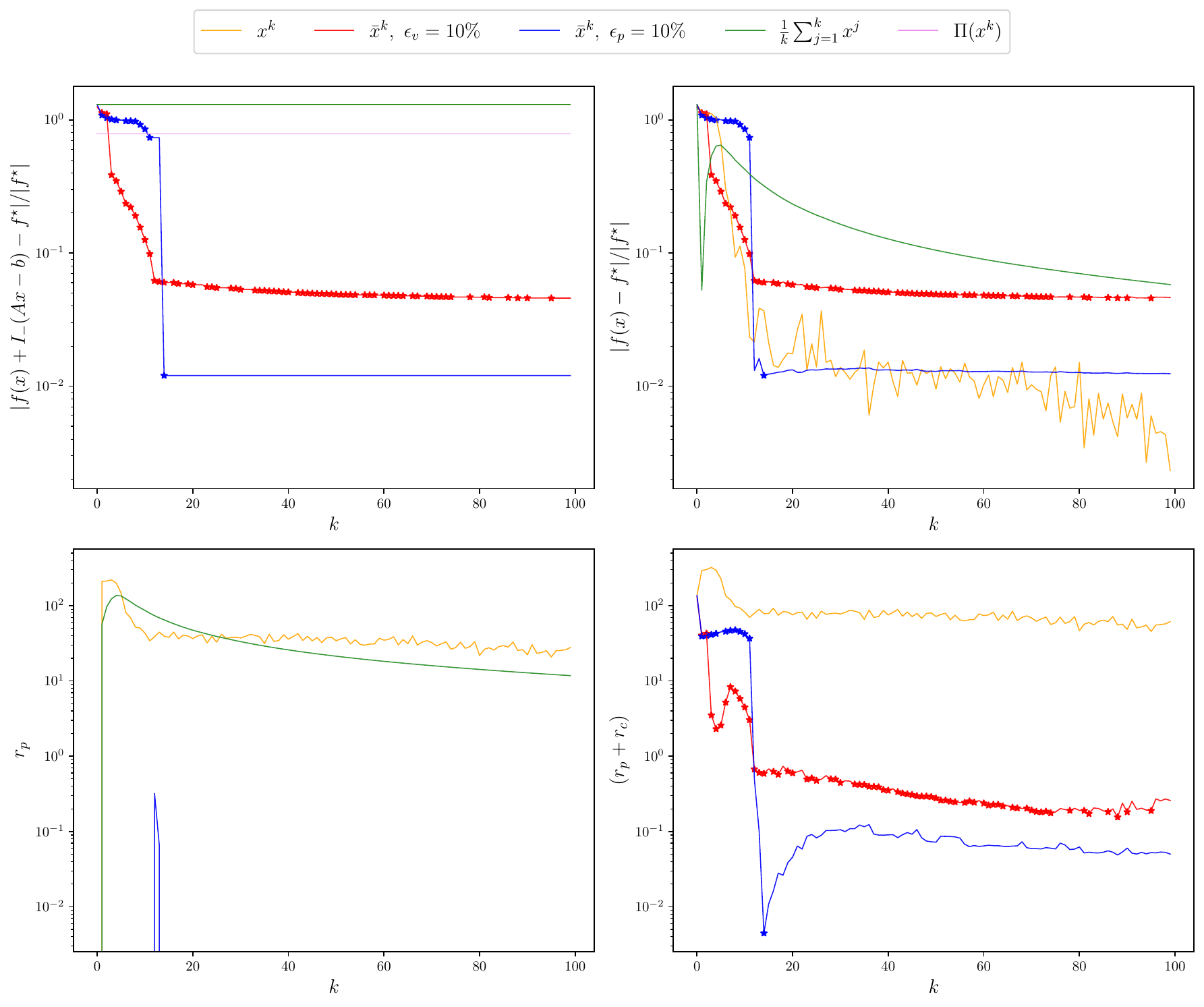}
    \caption{Best-to-date suboptimality of the most recent feasible point (top left),
    function value suboptimality (top right), relative primal violations (bottom left),
    relative residuals (bottom right)
    versus localization method iterations for
    convex relaxation of the assignment problem with history $H=3$.}
    \label{fig-cvx-ap-subopt-vp-h}
\end{center}
\end{figure}

\begin{figure}
\begin{center}
    \includegraphics[width=0.8\textwidth]{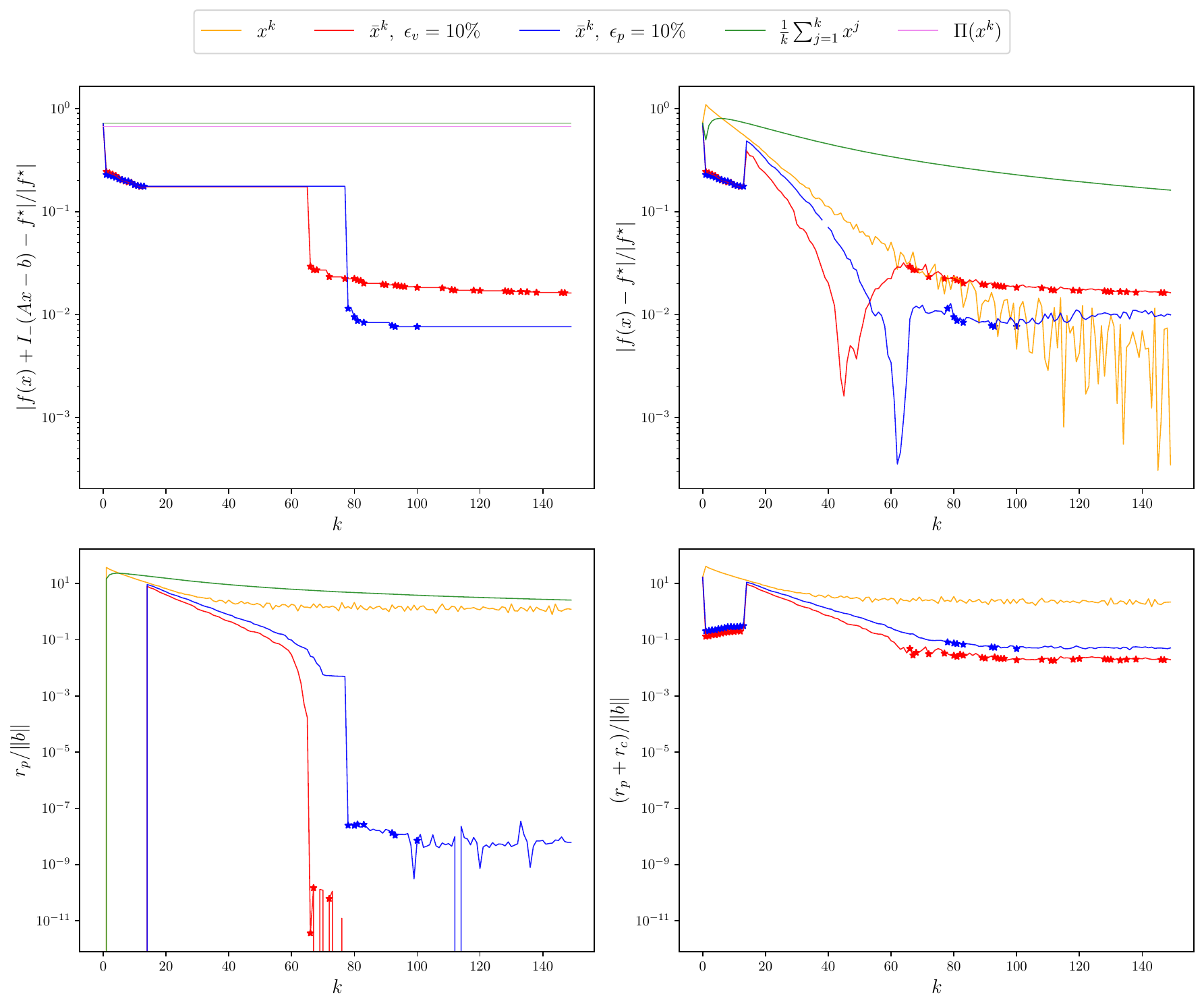}
    \caption{Best-to-date suboptimality of the most recent feasible point (top left),
    function value suboptimality (top right), relative primal violations (bottom left),
    relative residuals (bottom right)
    versus localization method iterations for
    the multi-commodity problem with history $H=5$.}
    \label{fig-mcf-subopt-vp-h}
\end{center}
\end{figure}

\begin{figure}
\begin{center}
    \includegraphics[width=0.8\textwidth]{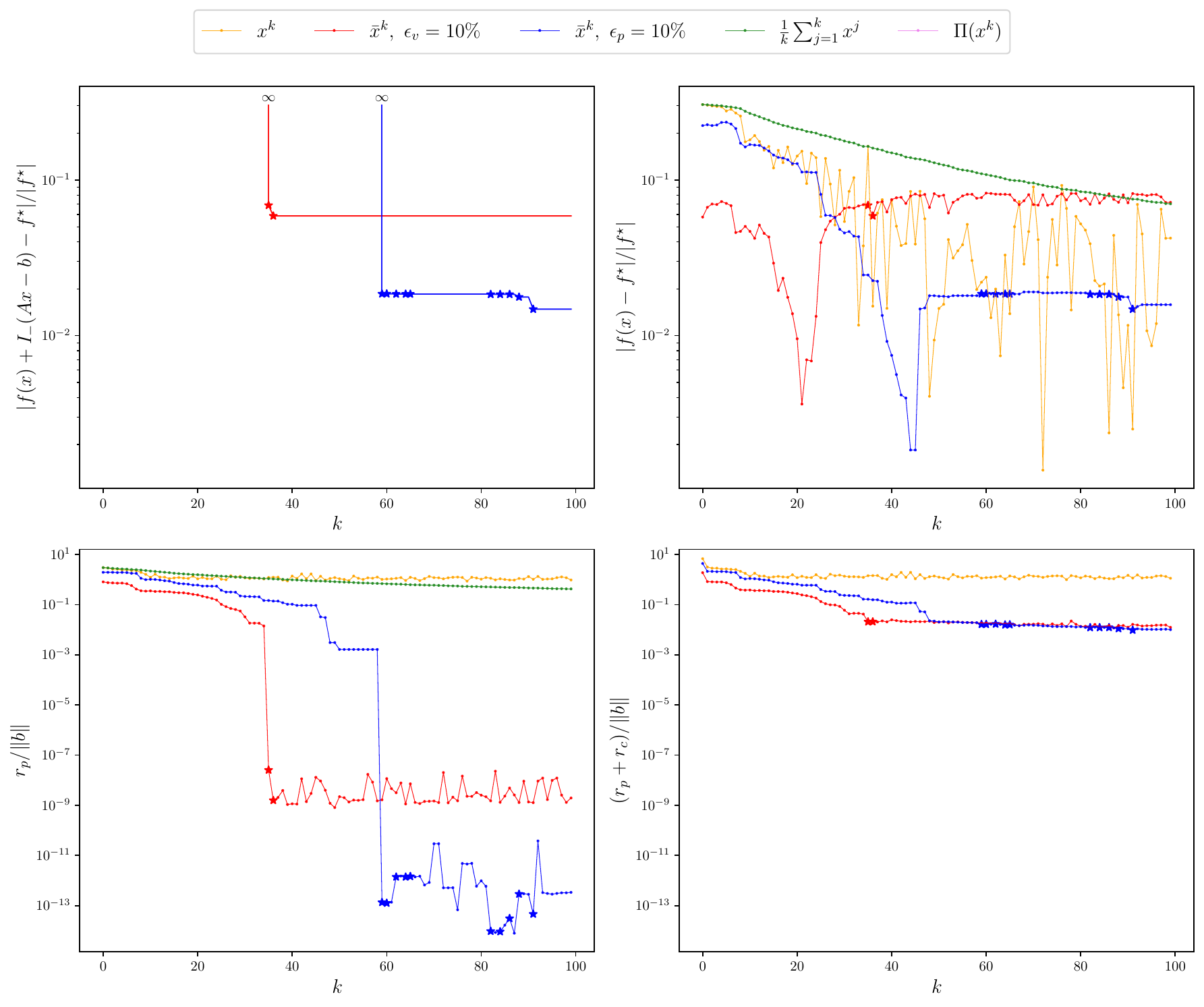}
    \caption{Best-to-date suboptimality of the most recent feasible point (top left),
    function value suboptimality (top right), relative primal violations (bottom left),
    relative residuals (bottom right)
    versus localization method iterations for
    the shipment problem with history $H=3$. }
    \label{fig-ot-subopt-vp-h}
\end{center}
\end{figure}

In this section, we test using historical data
for MARA primal recovery while solving problem \eqref{e-cvx-primal-recov-rp},
as explained in \S\ref{sec-MARA-extensions}. We used the localization method for price discovery.
Our experiments demonstrate that incorporating history $H$
accelerates the identification of primal feasible points, 
with a slight increase in suboptimality.
See Figures \ref{fig-ra-subopt-vp-h}, \ref{fig-cvx-ap-subopt-vp-h},
\ref{fig-mcf-subopt-vp-h}, and \ref{fig-ot-subopt-vp-h}.
For example, in the resource allocation problem  
$\epsilon_p=10\%$ primal recovery method with history $H=3$ finds a 
feasible point at iteration $37$ with suboptimality $0.25\%$,
instead of iteration $92$ with suboptimality $0.15\%$ without history ($H=1$),
see Figures~\ref{fig-ra-subopt-vp-h} and \ref{fig-ra-subopt-vp}, respectively.

Similarly, in the multi-commodity flow problem without history,
the $\epsilon_v=10\%$ MARA
does not recover the primal feasible point except the initial iteration, while
$\epsilon_p=10\%$ MARA recovers primal feasible point at iteration $137$
with suboptimality $0.14\%$.
In contrast, with history $H=5$, MARA with both dual methods recovers primal feasible 
points faster, at iterations $68$ and $78$,
with suboptimalities $2.7\%$ and $1.1\%$, respectively.
See Figures~\ref{fig-mcf-subopt-vp} and \ref{fig-mcf-subopt-vp-h}.

\section{Conclusions}
\label{sec-conclusions}

MARA works with any dual optimization algorithm, providing a flexible way to construct fast, feasible, and controllably-suboptimal solutions.
Numerical examples on several problem classes using both a price localization method and 
a dual subgradient method
show that MARA is able to generate feasible, near-optimal solutions in just a few tens of iterations of the underlying dual algorithm. While MARA increases total computation budgets significantly, because it operates independently of the dual algorithm using parallel calls to agents, it does not impact
the wall-clock time of the underlying dual algorithm. With the increasing scale of computing resource available in modern data centers, MARA's capability to leverage a large computing budget to quickly find feasible, near-optimal solutions for large-scale problems -- while providing hyperparameters
to control the speed and solution quality trade-off -- is valuable.

MARA is flexible, so further research will no doubt uncover other variations that can further improve its 
performance. Another useful extension is to non-convex agents, in particular agents with integrality constraints for which primal recovery is especially challenging. 
In this case, the MARA step of taking convex combinations of responses can be replaced by an integer program to select one response from each agent. (See Appendix \ref{sec-nonconvex}.) This an important practical extension on which we are actively working.

\clearpage
\bibliography{references}

\clearpage

\appendix

\section{Homogeneous ACCPM}\label{sec-haccpm}
The following properties are based on paper \cite{nesterov1999homogeneous}.
For any $z=(\bar \lambda, t)\in \reals_+^{m}\times\reals_{++}$, 
define $\lambda = \bar \lambda/t$ and
\[
\hat q(z) = (-g'(\lambda), g'(\lambda)^T \lambda) \in \reals^{m+1},
\qquad 
(-g)'(\lambda) = b - Ax(-A^T\lambda) \in \partial (-g)(\lambda).
\]
Then $q(z) = \hat q(z) / \|\hat q(z) \|_2$ is a homogeneous 
separation oracle with norm one.

Let the normal log-barrier be defined as
\[
F^{k+1}(z) = F^k(z) - \log (q(z_k)^T(z_k - z)).
\]
Define the following quantities
\[
\beta_{ik} = \frac{1}{q(z_i)^T(z_i - z_k)}>0, \quad i=1, \ldots, k,
\qquad S_k = \sum_{i=1}^k \beta_{ik},
\]
and the weighted average of slacks at $z$ as
\[
\mu_k(z) = \frac{1}{S_k} \sum_{i=1}^k \beta_{ik} q(z_i)^T(z_i - z).
\]
Then by Theorem 2.7 in \cite{nesterov1999homogeneous} we have
\[
\mu_k(z) \leq O({1}/{\sqrt{k}}).
\]

Also define the quantities
\[
\pi_{ik} = \beta_{ik}/\|\hat q(z_i)\|_2, \quad i=1, \ldots, k,
\qquad P_k = \sum_{i=1}^k \pi_{ik}.
\]
Using Lemma 4.3 in \cite{nesterov1999homogeneous}, we have
\[
\hat q(z)^T(z - \tilde z) = \tilde t (-g)'(\lambda)^T(\lambda - \tilde \lambda).
\]
Then for any $\lambda \in \reals^m_+$, and 
$z=(\lambda, 1)$, 
it follows
\[
\frac{1}{P_k} \sum_{i=1}^k \pi_{ik} (-g)'(\lambda_i)^T(\lambda_i - \lambda)
= \frac{1}{P_k} \sum_{i=1}^k \pi_{ik} \hat q(z_i)^T(z_i - z)
=\frac{S_k}{P_k} \mu_k(z) = O(1/\sqrt{k}),
\] 
see \cite{nesterov1999homogeneous}.
Define $\theta_{ik} = \pi_{ik}/P_k$ (clearly, $\sum_{i=1}^k \theta_{ik}=1$),  
this implies
\BEAS
O(1/\sqrt{k}) &\geq& \sum_{i=1}^k \theta_{ik} (-g)'(\lambda_i)^T(\lambda_i - \lambda^\star) \\
&\geq& g(\lambda^\star) - \sum_{i=1}^k \theta_{ik} g(\lambda_i) \\
&\geq& g(\lambda^\star) - g(\bar \lambda_k),
\EEAS 
where $\bar \lambda_k = \sum_{i=1}^k \theta_{ik} \lambda_i$.
In the above we used concavity of the dual function $g$ and Jensen inequality.
Combining the above with $g(\lambda^\star) \geq g(\bar \lambda_k)$, 
we get $g(\bar \lambda_k)\to g(\lambda^\star)$.

\section{Nonconvex agents}\label{sec-nonconvex}
MARA can be extended to the settings where the
problem~\eqref{e-prob} is nonconvex.
We assume that each agent has access to the
price-directed oracle~\eqref{e-distrib-oracle}.
Since the dual problem is always convex, any price discovery method
conforming to our access assumptions
(\eg, see \S\ref{s-price-local}) can be used directly to get optimal dual prices.
However, since the function $f$ is nonconvex, strong duality rarely holds.

\subsection{Lagrangian relaxation}
Since weak duality always holds, 
the dual problem is a Lagrangian relaxation of the nonconvex
primal problem.
Lagrangian relaxation has been studied since 1980s~\cite{shor1987quadratic}, 
for more results on duality see \cite[\S5]{boyd2004convex}
and \cite[\S5]{bertsekas1999nonlinear}.
The dual problem associated with problem~\eqref{e-prob} is given by
\BEQ\label{e-dual-problem}
\begin{array}{ll}
\mbox{maximize} & -f^*(u) - \lambda^T b\\ 
\mbox{subject to} & u = - A^T\lambda \\
                &   \lambda \geq 0. 
\end{array}
\EEQ
For $\lambda\geq 0$ and dual variable $\hat x$ associated with the equality constraints, 
the Lagrangian of~\eqref{e-dual-problem} is
\[
L(\lambda, u, \hat x) = -f^*(u) - \lambda^T b + \hat x^T(u + A^T\lambda).
\]
Thus the Lagrange dual problem associated with~\eqref{e-dual-problem} is 
\BEQ\label{e-second-dual}
\begin{array}{ll}
\mbox{minimize} & f^{**}(\hat x)\\ 
\mbox{subject to} & A\hat x \leq b.
\end{array}
\EEQ
By the conjugacy theorem~\cite{bertsekas2009convex}, the double conjugate 
function $f^{**}(\hat x)$ is a convex closure of function $f$, and since $f^*$ is a closed, convex, and proper function, we have
\[
(f^{**})^* = (f^*)^{**} = f^*.
\]
Therefore, assuming that there exists $\hat x \in \relint \dom f^{**}$ 
with $A\hat x\leq b$, then
by Slater's condition
strong duality holds between problems \eqref{e-second-dual} and \eqref{e-dual-problem}.

\subsection{Nonconvex primal recovery}
Due to the strong duality between problems \eqref{e-second-dual} 
and \eqref{e-dual-problem},
the optimality conditions for problem~\eqref{e-second-dual} 
are also~\eqref{e-opt-cond1}--\eqref{e-opt-cond4}.
Problem~\eqref{e-second-dual} is a convex relaxation of the nonconvex problem
\eqref{e-prob},
thus we generate a primal
variable reducing the primal residuals $r_p$ as follows.
Similar to \S\ref{sec-cvx-primal-recovery},
each agent $i=1,\ldots, K$ provides a list of $N_i$ $\epsilon$-suboptimal 
primal variables,
$\{z_{i_1}, \ldots, z_{N_i}\}$, 
associated with the current local price vector $y_i = - A_i^T\lambda$.
Since the functions $f_i$ are nonconvex, to produce primal variable 
$\bar x_i$ in the domain of $f_i$, 
we select a point from the set $\{z_{i_1}, \ldots, z_{N_i}\}$,
for all $i=1, \ldots, K$.
In particular, we solve the following mixed-integer linear program (MILP)
\BEQ\label{e-ncvx-primal-recov}
\begin{array}{ll}
\mbox{minimize} & r_p  \\
\mbox{subject to} & r_p = \ones^T (A\bar x-b)_+ \\
& \bar x_i = Z_i u_i, \quad i=1, \ldots, K \\
& \ones^T u_i = 1, \quad i=1, \ldots, K \\
& u_i \in \{0,1\}^{N_i}, \quad i=1, \ldots, K \\
& \bar x = (\bar x_1, \ldots, \bar x_K),
\end{array}
\EEQ
where $u_i$ is a variable,
and $Z_i$ is defined as in \eqref{e-zi-matrices},
for all $i=1, \ldots, K$.
The problem~\eqref{e-ncvx-primal-recov} can be readily solved using any MILP 
solver~\cite{achterberg2013mixed,santos2020mixed}.
Since this MILP may be computationally expensive, we next describe a simple heuristic for solving~\eqref{e-ncvx-primal-recov}.

\subsection{MILP heuristic}
Relax \eqref{e-ncvx-primal-recov} by replacing each nonconvex constraint $u_i\in\{0,1\}^{N_i}$ with its convex relaxation $u_i\in[0, 1]^{N_i}$
for all $i=1, \ldots, K$.
In the relaxed solution, interpret each component $u_{ij}$ as the success probability of the corresponding Bernoulli variable for all $j=1, \ldots, N_i$.
Next, we generate $S$ samples from these Bernoulli distributions and apply a greedy rounding method to each sample. 
Finally, we select the primal variable with the lowest residual 
as the approximate solution. We summarize the greedy rounding method:
\begin{enumerate}
	\item Begin with an initial set of fractional variables.   
    \item Successively iterate through each variable.
    \begin{enumerate} 
    \item For each variable, test the possible values $\{0,1\}$ while keeping 
    the other variables fixed.
    \item If any value results in a lower objective function, select that 
    value and proceed to the 
    next variable.
    \end{enumerate}
    \item Continue cycling until a local minimum is reached. 
\end{enumerate}

Note that step 2(a) can be carried out efficiently for
the objective \eqref{e-ncvx-primal-recov}. 


\end{document}

%% file: defs.tex
\newcommand{\ones}{\mathbf 1}

\newcommand{\reals}{{\mbox{\bf R}}}

\newcommand{\BEAS}{\begin{eqnarray*}}
\newcommand{\EEAS}{\end{eqnarray*}}
\newcommand{\BEA}{\begin{eqnarray}}
\newcommand{\EEA}{\end{eqnarray}}
\newcommand{\BEQ}{\begin{equation}}
\newcommand{\EEQ}{\end{equation}}
\newcommand{\BIT}{\begin{itemize}}
\newcommand{\EIT}{\end{itemize}}

\newcommand{\Tr}{\mathop{\bf Tr}}
\newcommand{\diag}{\mathop{\bf diag}}

\newcommand{\dom}{\mathop{\bf dom}} 
\newcommand{\relint}{\mathop{\bf relint}}

\newcommand{\eg}{{\it e.g.}}
\newcommand{\ie}{{\it i.e.}}

\newcounter{algorithmctr}[section]
\renewcommand{\thealgorithmctr}{\thesection.\arabic{algorithmctr}}
\newenvironment{algdesc}%
   {\refstepcounter{algorithmctr}
   \begin{list}{}{%
       \setlength{\rightmargin}{0\linewidth}%
       \setlength{\leftmargin}{.05\linewidth}}%
       \rmfamily\small
       \item[]{\setlength{\parskip}{0ex}\hrulefill\par%
        \nopagebreak{\bfseries\textsf{Algorithm \thealgorithmctr~}}}}%
   {{\setlength{\parskip}{-1ex}\nopagebreak\par\hrulefill} 
   \end{list}}

   {\refstepcounter{algorithmctr}
   \begin{list}{}{%
       \setlength{\rightmargin}{0\linewidth}%
       \setlength{\leftmargin}{.05\linewidth}}%
       \rmfamily\small
       \item[]{\setlength{\parskip}{0ex}\hrulefill\par%
        \nopagebreak{\bfseries\textsf{Algorithm}}}}%
   {{\setlength{\parskip}{-1ex}\nopagebreak\par\hrulefill} 
   \end{list}}